\documentclass[11pt]{article}
 \usepackage[colorlinks,citecolor=cyan,linkcolor=blue]{hyperref}
\usepackage{amsmath,amsfonts,amsbsy,amstext,amsopn,amssymb}
\usepackage{longtable}
\usepackage{theorem}
\usepackage{mathrsfs}
\usepackage{epsfig}
\usepackage{long table}
\usepackage{algorithm}
\usepackage{algpseudocode}
\usepackage{color}
\usepackage[numbers,sort&compress]{natbib}

\makeatletter \@addtoreset{equation}{section}

\makeatother

{\theorembodyfont{\slshape}
\newtheorem{theorem}{Theorem}[section]
\newtheorem{proposition}{Proposition}[section]
\newtheorem{lemma}{Lemma}[section]
\newtheorem{corollary}{Corollary}[section]

} {\theorembodyfont{\rmfamily}
\newtheorem{definition}{Definition}[section]

\newtheorem{remark}{Remark}[section]
\newtheorem{example}{Example}[section]
}
\def\proof{{\noindent\sc Proof. \quad}}

\def\eproof{{\mbox{}\hfill\qed}\medskip}

\newcommand{\R}{{\mathbb{R}}}

\newcommand{\cond}{{\mathscr    K }}

\newcommand{\rank}{{\rm{rank}}}

\newcommand{\diag}{{\rm{diag}}}

\newcommand{\vect}{{\sf{vec}}}

\def\Oh{{\mathcal O}}

\newcommand{\qed}{$\hfill{\Box}$} 

\def\Diag{{\rm  diag}}

\def\bsa{{\boldsymbol a}}
\def\bsb{{\boldsymbol b}}
\def\bsc{{\boldsymbol c}}
\def\bsd{{\boldsymbol d}}
\def\bse{{\boldsymbol e}}
\def\bsf{{\boldsymbol f}}
\def\bsh{{\boldsymbol h}}

\def\bs\omega{{\boldsymbol \omega}}
\def\bsp{{\boldsymbol p}}
\def\bsq{{\boldsymbol q}}

\def\bsg{{\boldsymbol g}}

\def\bsh{{\boldsymbol h}}

\def\bss{{\boldsymbol s}}
\def\bsv{{\boldsymbol v}}
\def\bsu{{\boldsymbol u}}
\def\bsw{{\boldsymbol w}}
\def\bst{{\boldsymbol t}}
\def\bsr{{\boldsymbol r}}

\def\bsl{{\boldsymbol l}}
\def\bsx{{\boldsymbol x}}
\def\bsy{{\boldsymbol y}}

\def\bsA{{\boldsymbol A}}
\def\bsB{{\boldsymbol B}}
\def\bsC{{\boldsymbol C}}
\def\bsD{{\boldsymbol D}}
\def\bsE{{\mathcal  E}}
\def\bsF{{\boldsymbol F}}
\def\bsR{{\boldsymbol R}}
\def\bsX{{\boldsymbol X}}

\def\bsI{{\boldsymbol I}}

\def\bsU{{\boldsymbol U}}

\def\bsL{{\boldsymbol L}}
\def\bsQ{{\mathcal   D}}

\def\bsS{{\boldsymbol S}}
\def\red{\textcolor[rgb]{0,0,0}}

\newcommand*{\affaddr}[1]{#1} 
\newcommand*{\affmark}[1][*]{\textsuperscript{#1}}

\begin{document}
\title{\bf Structured condition number for  multiple right-hand side linear systems with parameterized quasiseparable coefficient matrix}

\author{
Qingle Meng\affmark[1]
\thanks{Email:
qinglemeng@yahoo.com}
, \,
Huaian Diao\affmark[2]
\thanks{Corresponding author. Email: hadiao@nenu.edu.cn}
and   Qinghua Yu\affmark[2]
\thanks{Email: 17854103330@163.com}\\
\affaddr{
\affmark[1]School of Mathematical Sciences, Xiamen University, Xiamen 361005,\\P.R. China.}\\
\affmark[2]School of Mathematics and Statistics,
Northeast Normal University, \\
No. 5268 Renmin Street, Chang Chun 130024,
P.R. China.}

\date{}
\maketitle

\begin{quote}
{\small {\bf Abstract.} In this paper, we consider the structured perturbation analysis for multiple right-hand side linear systems with parameterized coefficient matrix. Especially, we present the explicit expressions for structured condition numbers for multiple right-hand sides linear systems with \{1;1\}-quasiseparable coefficient matrix in the quasiseparable and the Givens-vector representations. In addition, the comparisons of these two condition numbers between themselves, and with respect to unstructured condition number are investigated.\,\,Moreover, the effective structured condition number for multiple right-hand sides linear systems with \{1;1\}-quasiseparable coefficient matrix is proposed. The relationships between the effective structured condition number and structured condition numbers with respect to the quasiseparable and the Givens-vector representations are also studied.\,\,Numerical experiments show that there are situations in which the effective structured condition number can be much smaller than the unstructured ones.}
\end{quote}



\section{Introduction}
 The concepts of backward errors and condition numbers play an important role in numerical linear algebra \cite{Higham2002Book}. For example, to solve the linear system
\begin{equation}\label{eq:linear}
	\bsA \bsx =\bsb,
\end{equation}
the backward errors are used to measure the minimal  magnitude
 perturbations on the data $\bsA$ and $\bsb$ such that the computed solution is the exact solution to the perturbed linear system. On the other hand, condition numbers describe the worst-case sensitivity of the solution $\bsx $ of \eqref{eq:linear} to all possible small perturbations on both data $\bsA$ and data $\bsb$ \cite{Rice1966SIMNUA}; see the recent comprehensive monograph \cite{Cucker2013Book} for more details. The forward error of the computed solution can be bounded by the product of the backward error and condition number, which can tell us the loss of the procession of the numerical algorithm for solving the problem.

 \red{Many papers and books had been devoted to} the condition number theory for the linear system \eqref{eq:linear}; see \cite{GolubVanLoan2013Book,Higham2002Book} and their references therein. There are two types of condition numbers, called {\em normwise} and {\em componentwise} condition numbers; see \cite[Chap.7]{Higham2002Book}. The normwise condition numbers measure the input and output error by means of norms, while the componentwise condition numbers use componentwise perturbations to measure the error of the input data.
  We should point out that when the data is sparse and badly-scaled, it is more suitable to adopt the componentwise condition number to define the conditioning of the problem, because the normwise condition numbers \red{ignore} the sparsity and scaling of the data \red{and} may overestimate the exact conditioning of the problem. For structured linear systems, it is \red{reasonable} to investigate structured perturbations on the input data, because structure-preserving algorithms that preserve the underlying matrix structure can enhance the accuracy and efficiency of solving linear systems. The structured condition numbers for structured linear systems can be found in \cite{Dopico2016,Higham1992Str,Rump03a,Rump03b}.

 In this paper, we  focus on  the linear system with multiple right-hand sides
\begin{equation}\label{eq:mul}
	\bsA \bsX = \bsB,
\end{equation}
where $\bsA$ $\in \R^{n\times n}$ is nonsingular and $\bsB$ $\in \R^{n \times m}$. The above equation is a generalization of \red{linear system} \eqref{eq:linear}.
\red{The multiple right-hand side linear system \eqref{eq:mul} is indeed a special case of the linear matrix equation--Sylvester equation
\begin{equation}\label{eq:syl}
	\bsA\bsX+\bsX\bsC=\bsB
\end{equation}
by taking $\bsC=0$. Matrix equations including Sylvester equations and algebraic Riccati equations  have many applications in problems of control \cite{Antoulas}, discretization of PDE \cite{Palitta}, block-diagonalization \cite{GolubVanLoan2013Book}, and many others. For recent numerical algorithm developments for solving Sylvester equations and algebraic Riccati equations, we refer to papers  \cite{Bai2011,Bai2006,GaoAndBai2011}. Especially,  the multiple right-hand side linear system \eqref{eq:mul} arises naturally in many  applications such as Quantum Chromo Dynamics \cite{SakuraiAndTadano2010}, dynamics of structures \cite{CloughAndPenzien1975}, quasi-Newton methods for solving nonlinear equations with multiple secant equations \cite{Higham1992Componentwise}, computing the lengths of nucleon-nucleon scattering \cite{StathopoulosAndOrginos2010}, wave propagation phenomena \cite{VirieuxAndOperto2009}.} For perturbation analysis for (generalized) Sylvester equation, *-Sylvester equation and algebraic Riccati equation, we refer to papers \cite{diao12,diao13,diao17,diao17*Sy,Higham93}.

The low-rank structured matrix has been studied extensively in numerical linear algebra and \red{has} many applications; see the recent books \cite{EidelmanETC2014Linear,EidelmanETC2014Eig} and references therein. The large submatrices of the low-rank structured matrix have ranks much smaller than the the size of the matrix. Based on the above observation, the $n$-by-$n$ low-rank structured matrix can be represented by different sets of $\Oh(n)$ parameters, which is named as {\em representations} \cite[Chap. 2]{VandebrilETC2008Linear}. This kind of representations can help us to develop {\em fast algorithms}, which directly operate on the parameters of the representation, with the computational costs of $\Oh(n)$ for solving linear systems or $\Oh(n^2)$ for solving eigenvalue problems. \red{{\em Semiseparable matrices} \cite{eidelman1999on} are a special category of low-rank structured matrices, which appear in several types of applications, e.g. the field of integral equations \cite{Geonzalesa1997,Petersen2004}, boundary value problems \cite{Geonzalesa1997,Greengard1991,Lee1997,HP1992}, in the theory of Gauss-Markov processes time-varying linear systems \cite{Dewilde1998,Gohberg1984}, acoustic and electromagnetic scattering theory \cite{Colton1998} and rational interpolation \cite{Barel2002}. Furthermore, it was shown that in \cite{vandebril2005a}  semiseparable matrices have other equivalent representations which are named as  {\em quasiseparable matrices}. Faster solver for the low-rank structured Sylvester and Lyapunov equations has been proposed in  \cite{Palitta2018}. It was shown that for Sylvester equations \eqref{eq:syl},   when $\bsA$, $\bsB$ and $\bsC$ are quasiseparable, the solution $\bsX$ is numerically quasiseparable. Therefore, it is interesting to consider the conditioning of the solution $\bsX$ with respect to the representations of the low-rank structured coefficient matrices.}

\red{Recently, structured componentwise condition numbers for low-rank structured matrices have been introduced for eigenvalue problems \cite{Dopico2015Structured}, linear systems \cite{Dopico2016}, and generalized eigenvalue problems \cite{diaomeng} with parameterized quasiseparable matrices \cite{vandebril2005a}. In this paper, we consider the structured componentwise condition numbers of the multiple right-hand side  linear system \eqref{eq:mul} when the coefficient matrix $\bsA$ is a quasiseparable matrix. Especially, we focus on the case that $\bsA$ is a \{1;1\}-quasiseparable matrix. The class of \{1;1\}-quasiseparable matrices incorporates the class of semiseparable matrices of semiseparability rank, tridiagonal matrices and unitary Hessenberg matrices etc.;  see \cite{Geonzalesa1997,Petersen2004} for details. In this paper, we mainly focus on two most important representations for \{1;1\}-quasiseparable coefficient matrix $A$. One is the general {\em quasiseparable representation} \cite{VandebrilETC2008Linear}, which is non unique, and another is the essentially unique {\em Givens-vector representation} \cite{vandebril2005a}, which is introduced to improve the numerical stability of fast matrix computations involving quasiseparable matrices.} The explicit expression of structured componentwise condition numbers for \red{the multiple right-hand side  linear system} \eqref{eq:mul} with respect to these two representations is derived. \red{We should point out that recent works \cite{Amestoy2015f,Amestoy2017f} have been done by exploiting the sparsity of of the multiple right-hand sides to reduce the computational cost of the direct solver for \eqref{eq:mul}. The linear system \eqref{eq:mul} with sparse multiple right-hand sides arises in several real-life applications, for example, 3D frequency-domain full waveform inversion and 3D
controlled-source electromagnetic inversion \cite{Mary2017}. Therefore, it is suitable to consider condition numbers for the solution to  \eqref{eq:mul} with respect to the sparsity of $\bsB$.}  In this paper, we propose the effective structured componentwise condition number for the solution of \eqref{eq:mul}, which can be used to measure the structured conditioning of \eqref{eq:mul} with respect to the structured perturbations on the coefficient matrix and sparse multiple right-hand sides. Numerical examples show that there are bigger differences between the structured effective condition number and the unstructured ones, which means that \eqref{eq:mul} is less sensitive to the structured perturbations on $\bsA$ and $\bsB$. The above observation can help us to understand why structure-preserving algorithms are important \red{in solving the multiple right-hand side linear system} \eqref{eq:mul} in addition to reducing the computational complexity, since the forward error of the algorithm can be reduced significantly with respect to the structured componentwise perturbations on the coefficient matrix $\bsA$ and sparse multiple right-hand sides $\bsB$.




  This paper is organized as follows.\,\,The basic results of \{1;1\}-quasiseparable matrices are reviewed in Section \ref{preli:basic-result}. Different types of condition numbers for multiple right-hand sides linear systems with general parameterized coefficient matrices are investigated in Section \ref{fundemental def}, which can be used to derive explicit expressions for the condition number of multiple right-hand side linear systems with \{1;1\}-quasiseparable coefficient matrix with respect to the quasiseparable representation \cite{eidelman1999on} and the Givens-vector representation via tangent \cite{vandebril2005a,Dopico2015Structured,Dopico2016} in Section \ref{section:quasi-represetation}. We study relationships between different structured and unstructured condition numbers for multiple right-hand side linear systems with parameterized coefficient matrix in Section \ref{section:6}. Numerical experiments are implemented for random generated parameters defining $\bsA$ and $\bsB$ in Section \ref{section:8}. Concluding remarks are drawn and future research topics are pointed out in Section \ref{sec:co}.

 {\bf Notations.} In this paper, we adopt the following notations. For any two conformal matrices $\bsA $ and $\bsB$ \red{of size $n\times n$}, $\bsA \leq  \bsB$ should be understood componentwise. The notation \red{$0\leq \bsA$ implies  that $\bsA$ is nonnegative, and the similar notation is adopted for vectors.} $\bsB/\bsA:=(\bsb_{i,j}/\bsa_{i,j} )$ where $\bsa_{i,j}$ is $(i,j)$-th entry of $\bsA$ and if $\bsa_{i,j}=0$ the corresponding $(i,j)$-th entry of $\bsB/\bsA$ and  $\bsB$ should be zero,  $\bsA^\top $ is the transpose matrix of $\bsA$, $|\bsA|$ is obtained by taking absolute values operation on each entry of $\bsA $, $\vect(A)$ is a column vector formed by  the columns of $\bsA  $ one by one, $\bsA(:,i)$ and $\bsA(j,:)$ are $i$-th column and $j$-th row of $\bsA$ respectively, $\bsA(i_1:i_2,j_1:j_2)$ is a submatrix of $\bsA \in \R^{n\times n}$ consisting of rows $i_1$ up to and including $i_2$ and columns $j_1$ up to and including $j_2$ of $\bsA$ with $1\leq i_1\leq i_2\leq n$ and $1\leq j_1\leq j_2 \leq n$, $\Diag(\bsa) \in  \R^{p \times p}$ is a diagonal matrix with the diagonal entries being the corresponding entres of $\bsa \in \R^{p \times 1}$, $\|\bsx\|_\infty$ is the $\infty$-norm of a vector $\bsx$, $\|\bsA\|_F$ is the Frobenius norm of $\bsA$, $\| \bsA\|_{\max }=\max\limits_{i,j} |\bsa_{i,j}|$. \red{The inequality $\mathbf{0}\neq \bsX$ means that  $0\neq \bsX(i,j)$ for all $i,j$.  }

%

\section{Preliminaries}\label{preli:basic-result}

In this section, we will give a brief review about \{1;1\}-quasiseparable matrices, which \red{are} a particular class of quasiseparable matrices. Especially, the general quasiseparable representation and Givens-vector representation for \{1;1\}-quasiseparable matrices are introduced. Moreover, the explicit derivatives expressions of entries of the \{1;1\}-quasiseparable matrix $\bsA$ with respect to its parameters in the general quasiseparable representation and Givens-vector representation of $\bsA$ are reviewed, which will \red{help} us to derive the componentwise condition numbers for multiple right-hand side linear systems with \{1;1\}-quasiseparable coefficient matrix $\bsA$ with respect to the parameters defining $\bsA$.

Quasiseparable matrices \red{introduced firstly in \cite{eidelman1999on}} are particular cases of low-rank structured matrices. The \{1;1\}-quasiseparable matrix is a special subclass of quasiseparable matrices, which satisfies that $\max_{i} \rank( \bsA(i+1:n,1:i ))= 1$ and $\max_{i} \rank (\bsA(1:i,i+1:n)) = 1$. Moreover, the maximum rank of the lower and upper submatrix  equal to $1$.
In this paper we adopt the representation of an $n$-by-$n $ \{1;1\}-quasiseparable matrices with $\Oh(n)$ parameters instead of its $n^{2}$ entries; see more details in \cite{Dopico2015Structured,Dopico2016}.


\begin{definition}\label{th:{1;1}}

A matrix $\bsA$ $\in \R^{n\times n}$ is a \{1;1\}-quasiseparable matrix if and only if it can be parameterized in terms of the following set of $7n-8$ real parameters,
  \begin{align}\label{eq:qs general}
   \Omega_{QS}=\Big(\{\bsp_{i}\}_{i=2}^{n},\{\bsa_{i}\}_{i=2}^{n-1},\{\bsq_{i}\}_{i=1}^{n-1},\{\bsd_{i}\}_{i=1}^{n},\{\bsg_{i}\}_{i=1}^{n-1},\{\bsb_{i}\}_{i=2}^{n-1},\{\bsh_{i}\}_{i=2}^{n}\Big),
  \end{align}
  as follows,
$$\bsA=      \begin{bmatrix}
          \bsd_{1} &\bsg_{1}\bsh_{2}&\bsg_{1}\bsb_{2}\bsh_{3} & \cdots &\bsg_{1}\bsb_{2}\cdots \bsb_{n-1}\bsh_{n}\\
          \bsp_{2}\bsq_{1} &\bsd_{2}&\bsg_{2}\bsh_{3} & \cdots &\bsg_{2}\bsb_{3}\cdots \bsb_{n-1}\bsh_{n}\\
          \bsp_{3}\bsa_{2}\bsq_{1} &\bsp_{3}\bsq_{2}&\bsd_{3} & \cdots &\bsg_{3}\bsb_{4}\cdots \bsb_{n-1}\bsh_{n}\\
          \bsp_{4}\bsa_{3}\bsa_{2}\bsq_{1} &\bsp_{4}\bsa_{3}\bsq_{2}&\bsp_{4}\bsq_{3} & \cdots &\bsg_{4}\bsb_{5}\cdots \bsb_{n-1}\bsh_{n}\\
          \vdots &\vdots &\vdots &\ddots &\vdots\\
          \bsp_{n}\bsa_{n-1}\cdots \bsa_{2}\bsq_{1} &\bsp_{n}\bsa_{n-1}\cdots \bsa_{3}\bsq_{2}&\bsp_{n}\bsa_{n-1}\cdots \bsa_{4}\bsq_{3} & \cdots &\bsd_{n}\\
         \end{bmatrix}.
 $$
 \end{definition}

In the following lemma, we will review the explicit derivative expressions of entries of the \{1;1\}-quasiseparable matrix  $\bsA$ with respect to the parameters $\Omega_{QS}$ given in Definition \ref{th:{1;1}}.

\begin{lemma}\label{lemm:proper-all-param}\cite[Theorem 4.4]{Dopico2015Structured}
Let $\bsA$ $\in \R^{n \times n} $ be a \{1;1\}-quasiseparable matrix and $\bsA=\bsA_{\bsL}+\bsA_{\bsD}+\bsA_{\bsU}$, with $\bsA_{\bsL}$ strictly lower triangular,
$\bsA_{\bsD}$ diagonal, and $\bsA_{\bsU}$ strictly upper triangular. Let $\Omega_{QS}$ be a quasiseparable representation of $\bsA$ given in Definition \ref{th:{1;1}}.
Then the entries of $\bsA$ are differentiable functions of the parameters in $\Omega_{QS}$ and
\begin{align*}
\left\{ \begin{array}{ll}
	 {\rm (1)} \quad \dfrac{\partial\bsA}{\partial\red{\bsd_{i}}} = e_{i}e^{\top}_{\red{i}}, &  i=1,2,\ldots,n,\\[10pt]
       {\rm (2)} \quad \bsp_{i}\dfrac{\partial\bsA}{\partial\bsp_{i}}=e_{i}\bsA_{\bsL}(i,:), & i=2,\ldots,n,\\[10pt]
       {\rm (3)}\quad  \bsa_{i}\dfrac{\partial\bsA}{\partial\bsa_{i}}=\begin{bmatrix}     0&0\\\bsA(i+1:n,1:i-1)&0 \end{bmatrix},  &i=2,\ldots,n-1,\\
      {\rm (4)} \quad \bsq_{i}\dfrac{\partial\bsA}{\partial\bsq_{i}}=\bsA_{\bsL}(:,i)e^{\top}_{i}, & i=1,2,\ldots,n-1,\\[10pt]
        {\rm (5)} \quad \bsg_{i}\dfrac{\partial\bsA}{\partial\bsg_{i}}=e_{i}\bsA_{\bsU}(i,:),  & i=1,2,\ldots,n-1,\\[10pt]
        {\rm (6)} \quad \bsb_{i}\dfrac{\partial\bsA}{\partial\bsb_{i}}=\begin{bmatrix}
                                                                  0&\bsA(1:i-1,i+1:n)\\0&0
                                                                  \end{bmatrix}, & i=2,\ldots,n-1,\\[10pt]
       {\rm (7)} \quad  \bsh_{i}\dfrac{\partial\bsA}{\partial\bsh_{i}}=\bsA_{\bsU}(:,i)e^{\top}_{i}, &  i=2,\ldots,n\red{,}
\end{array}
 \right.
       \end{align*}
       where $e_i$ is the $i$-th column of the $n$-by-$n$ identity matrix.
\end{lemma}

Another important representation for the \{1;1\}-quasiseparable matrices, which is named as {\em Givens-vectors representation}, was introduced in \cite{vandebril2005a,Dopico2015Structured,Dopico2016}. \red{Givens-vectors representation can improve the numerical stability of fast computations involving quasiseparable matrices.} In the following definition we will recall this type representation.




\begin{definition}\label{th:givens repre}

A matrix $\bsA$ $\in \R^{n\times n}$ is a \{1;1\}-quasiseparable matrix, if and only if it can be parameterized in terms of the following set of real parameters.
\begin{description}
\item[{$\bullet$}] $\{\bsc_{i},\bss_{i}\}$, where $(\bsc_{i},\bss_{i})$ is a pair of consine-sine with $\bsc_{i}^{2}+\bss_{i}^{2}=1$ for every $i\in\,\, \{2,3,\cdots,n-1\}$,
\item[{$\bullet$}] $\{\bsv_{i}\}_{i=1}^{n-1}$,$\{\bsd_{i}\}_{i=1}^{n}$,$\{\bsw_{i}\}_{i=1}^{n-1}$, where all of them are independent real parameters,
\item[{$\bullet$}] $\{\bsr_{i},\bst_{i}\}_{i=2}^{n-1}$, where $(\bsr_{i},\bst_{i})$ is a pair of cosine-sine with $\bsr_{i}^{2}+\bst_{i}^{2}=1$
                    for every $i\in \{2,3,\cdots,n-1\}$,
\end{description}
as follows:
$$
\bsA= \begin{bmatrix}
          \bsd_{1} &\bsw_{1}\bsr_{2} & \cdots &\bsw_{1}\bst_{2}\cdots \bst_{n-2}\bsr_{n-1}&\bsw_{1}\bst_{2}\cdots \bst_{n-1}\\
          \bsc_{2}\bsv_{1} &\bsd_{2}& \cdots &\bsw_{2}\bst_{3}\cdots \bst_{n-2}\bsr_{n-1}&\bsw_{2}\bst_{3}\cdots \bst_{n-1}\\
          \bsc_{3}\bss_{2}\bsv_{1} &\bsc_{3}\bsv_{2}& \cdots &\bsw_{3}\bst_{4}\cdots \red{\bst_{n-2}}\bsr_{n-1}&\bsw_{3}\bst_{4}\cdots \bst_{n-1}\\
          \vdots &\vdots &\vdots &\ddots &\vdots\\
          \bsc_{n-1}\bss_{n-2}\cdots \bss_{2}\bsv_{1} &\bsc_{n-1}\bss_{n-2}\cdots \bss_{3}\bsv_{2}& \cdots&\bsd_{n-1} &\bsw_{n-1}\\
          \bss_{n-1}\bss_{n-2}\cdots \bss_{2}\bsv_{1} &\bss_{n-1}\bss_{n-2}\cdots \bss_{3}\bsv_{2}& \cdots &\bsv_{n-1}&\bsd_{n}\\
         \end{bmatrix}.
 $$
This representation is denoted by $\Omega_{QS}^{GV}$, i.e., $\Omega_{QS}^{GV}:=(\{\bsc_{i},\bss_{i}\}^{n-1}_{\red{i=2}},\{\bsv_{i}\}^{n-1}_{i=1},\{\bsd_{i}\}_{i=1}^{n},\{\bsw_{i}\}^{n-1}_{i=1},\\ \{\bsr_{i},\bst_{i}\}^{n-1}_{i=2})$.

\end{definition}


 Because of the \red{dependency of $c_{i}$ and $s_{i}$}, which is also true for $\{r_{i},t_{i}\}$, there are correlated parameters in the Givens-vector representation. In this paper, we will focus on arbitrary componentwise perturbations of $\Omega_{QS}^{GV}$, which obviously destroy the cosine-sine pairs, thus it is more reasonable to restrict perturbations that preserve \red{such relationships of the cosine-sine pairs}. In \cite{Dopico2015Structured}, additional parameters of these pairs $\{c_{i},s_{i}\}$ and $\{r_{i},t_{i}\}$ are introduced by using their corresponding tangents, which are called as the {\em Givens-vector representation via
tangent} for the \{1;1\}-quasiseparable matrix.

\begin{definition}\label{def:givens }
For any Givens-vector representation $\Omega^{GV}_{QS}$ given in Definition \ref{th:givens repre}
of \{1;1\}-quasiseparable matrix  $\bsA$$\in \R^{n\times n}$, we define the {\em Givens-vector representation via tangent} as
\begin{align}\label{eq:GV rep}
\Omega_{GV}:&=\Big(\{\bsl_{i}\}^{n-1}_{i=2},\{\bsv_{i}\}^{n-1}_{\red{i=1}},\{\bsd_{i}\}^{n}_{i=1},\{\bsw_{i}\}^{n-1}_{i=1},\{\bsu_{i}\}^{n-1}_{i=2},\Big),
\end{align}
where
\begin{align}
\bsc_{i}&= \frac{1}{\sqrt{1+\bsl^{2}_{i}}},\quad \bss_{i}=\frac{\bsl_{i}}{\sqrt{1+\bsl^{2}_{i}}}, \quad \bsr_{i}=\frac{1}{\sqrt{1+\bsu^{2}_{i}}},\quad \bst_{i}=\frac{\bsu_{i}}{\sqrt{1+\bsu^{2}_{i}}}, \mbox{for } i=2,3,\ldots,n-1.\notag
\end{align}
\end{definition}
%

Lemma \ref{lemm:derivati-partial-Givens} gives explicit derivative expressions of entries of the \{1;1\}-quasiseparable matrix $\bsA$  with respect to the tangent-Givens-vector representation, which appeared in \cite{Dopico2015Structured}.


\begin{lemma}\label{lemm:derivati-partial-Givens}
Let $\bsA$ $\in \R^{n\times n }$ be a \{1;1\}-quasiseparable matrix and let $\Omega_{QS}$ be the tangent-Givens- vector representation of $\bsA$, where $\Omega_{GV}=\big( \{\bsl_{i}\}^{n-1}_{i=2}, \{\bsv_{i}\}^{\red{n-1}}_{i=1}, \{\bsd_{i}\}^{n-1}_{i=1}, \{\red{\bsw_{i}}\}^{n-1}_{i=1}, \{\bsu_{i}\}^{n-1}_{i=2} \big) $. Then the entries of $\bsA$ are differentiable functions of the parameters in $\Omega_{GV}$ and
\begin{align*}
\left\{
\begin{array}{ll}
	{\rm (1)} \quad \bsl_{i} \dfrac{\partial \bsA}{\partial\bsl_{i}} = \begin{bmatrix} 0&0\\-\bss^{2}_{i}\bsA(i,1:i-1)&0\\\bsc^{2}_{i}\bsA(i+1:n,1:i-1)&0\\
 \end{bmatrix}, & i=2,3,\ldots,n-1, \\[10pt]
{\rm (2)} \quad \bsu_{i} \dfrac{\partial \bsA}{\partial\bsu_{i}} = \begin{bmatrix} 0&-\bst^{2}_{i}\bsA(1:i-1,i+1:n)&\red{\bsr^2_{i}\bsA(1:i-1,i+1:n)}\\0&0&0\\
 \end{bmatrix},  & i=2,3,\ldots,n-1.
\end{array}
\right.
 \end{align*}
 \end{lemma}


\begin{remark}\label{re:dve}
Recalling that
$
\{\bsd_{i}\}^{n}_{i=1},\{\bsv_{i}\}^{n-1}_{\red{i=1}},\{\bsw_{i}\}^{n-1}_{i=1}
$
in the tangent-Givens-vector representation of $\bsA $ can be viewed as
$
\{\bsd_{i}\}^{n}_{i=1},\{\bsq_{i}\}^{n-1}_{i=2},\{\bsg_{i}\}^{n-1}_{i=1}
$
in a quasiseparable representation of $\bsA $, we have a similar expressions of the partial derivatives with respect to the parameters in $\{\bsd_{i}\}^{n}_{i=1},\{\bsv_{i}\}^{n-1}_{\red{i=1}},\{\red{\bsw_{i}}\}^{n-1}_{i=1}$ as 1), 4) and 5) in Lemma \ref{lemm:proper-all-param}.
\end{remark}

\section{Condition numbers for multiple right-hand side  linear systems with general parameterized coefficient matrix }\label{fundemental def}

In this section, we will first review unstructured componentwise condition number for \eqref{eq:mul}. \red{A general parameterized representation for $\bsA$ and $\bsB$ in \eqref{eq:mul} is introduced in Definition \ref{def:com-paramter-condi-num}, where we assume that there are some common parameters defining $\bsA$ and $\bsB$ simitounously. In Theorem \ref{th:explict-express-com-param-con-num1}, we will introduce the structured componentwise condition number with respect to parameters introduced in Definition \ref{def:com-paramter-condi-num}. We apply the result of Theorem \ref{th:explict-express-com-param-con-num1} to the case that $\bsB$ is a general dense matrix in Theorem \ref{th:explict-express-com-param-con-num2}, i.e., the parameter representation of $\bsB$ are just its entries. The multiple right hand sides $\bsB$ is usually sparse in some practical applications, for example, 3D frequency-domain full waveform inversion and 3D controlled-source electromagnetic inversion \cite{Mary2017}. Therefore, when $\bsB$ is sparse, we derive the structured componentwise condition number of the solution to \ref{eq:mul} with respect to the sparse pattern of $\bsB$ and parameterized representation of $\bsA$ in Theorem \ref{th:explict-express-com-param-con-num}. When  $\bsA$ is a general unstructured matrix and $\bsB$ is sparse, the componentwise condition number of \eqref{eq:mul} with respect to  the parameter representation of $\bsB$  and entries of $\bsA$ is given in Corollary \ref{co:3.1}.  }

In Definition \ref{def:com-con_number} we first review the componentwise condition number $\cond_{\bsE,\bsF}(\bsA,\bsB)$ for \eqref{eq:mul}, which was introduced in \cite[Definition 2.2]{zhang2007condition}. Furthermore, the explicit expression of $\cond_{\bsE,\bsF}(\bsA,\bsB)$ appeared in  \cite[Theorem 2.3]{zhang2007condition}. 

\begin{definition}\label{def:com-con_number}\cite{zhang2007condition}
Let $\bsA \bsX =\bsB$, where $\bsA$ $\in \R^{n\times n}$ is nonsingular, $\bsB$ $\in \R^{n \times m}$ and $\mathbf{0}\neq \bsX$ $\in \R^{n \times m}$. Then, for $\mathbf{0}\leq \bsE$ $ \in\R^{n \times n}$ and $\mathbf{0}\leq \bsF$ $ \in\R^{n \times m}$, we define the componentwise  condition number as
\begin{align}
\cond_{\bsE,\bsF}(\bsA,\bsB)&:=\lim_{\eta\rightarrow0}\sup\Bigg \{\frac{\|\delta\bsX\|_{\max}}{\eta \|\bsX\|_{\max}}:(\bsA+\delta\bsA)(\bsX+\delta\bsX)=\bsB+\delta\bsB,\notag \\
&\hspace{4.5cm} |\delta\bsA|\leq\eta\,\bsE, |\delta\bsB|\leq\eta\,\bsF                                                                                 \Bigg \} \notag \\
	&=\frac{\Big\|\left|\bsA^{-1}\right|\, \bsE\, |\bsX|+\left|\bsA^{-1}\right |\, \bsF\Big\|_{\max}}{\|\bsX\|_{\max}}. \label{th:the expression-compon}
\end{align}
\end{definition}
The weight matrices $\bsE$ and $\bsF$ are flexible, for example if we take $\bsF=\bf 0$ there is no perturbations on the multiple right-hand sides.  A natural choice of $\bsE$ and $\bsF$ is $|\bsA|$ and $|\bsB|$, respectively. \red{Furthermore, if $\bsE=|A|$ and $\bsF=|B|$, from \eqref{th:the expression-compon} we have
\begin{align}
\cond_{|\bsA|,|\bsB|}(\bsA,\bsB)
	&=\frac{\Big\|\left|\bsA^{-1}\right|\, |\bsA|\, |\bsX|+\left|\bsA^{-1}\right |\, |\bsB|\Big\|_{\max}}{\|\bsX\|_{\max}}. \label{th:1the expression-compon}
\end{align}
}




\begin{remark}\label{remark:Ax=b-expression-Compon}
For the case \red{$\bsB=\bsb$\,(i.e., $m=1$)}\,is a vector,  \eqref{th:the expression-compon} degenerates into the relative componentwise  condition number for the linear system $\bsA \bsx =\bsb$ \red{defined by}:
\begin{align*}
\cond_{\bsE,\bsy  }(\bsA,\bsb)&:=\lim_{\eta\rightarrow0}\sup\Bigg \{\frac{\|\delta
\bsx \|_{\infty } } {\eta \|\bsx\|_{\infty}}
:(\bsA+\delta\bsA)(\bsx+\delta\bsx)=\bsb+\delta \bsb,  |\delta\bsA|\leq\eta\bsE, |\delta\bsb|\leq\eta \, \bsy                                                                                   \Bigg \},\nonumber \\
	&=\frac{\Big\|  |\bsA^{-1}|\bsE|\bsx|+|\bsA^{-1}|\bsy  \Big\|_{\infty}}{\|\bsx\|_{\infty}}.
\end{align*}
 See \cite[Chap.7]{Higham2002Book} for more details.  Here $\bsy \geq 0$.
\end{remark}

%
%

Because many interesting classes of matrices can be represented by sets of parameters instead of their entries, it is interesting to consider the componentwise perturbations for these parameters, which lead to the following componentwise condition number \eqref{def:cond_parame}  for \eqref{eq:mul}, \red{where both the coefficient matrix $\bsA$ and multiple right-hand sides $\bsB$  have parametric representations. In Definition \ref{def:com-paramter-condi-num},  we consider a general situation that $\bsA $ and $\bsB$ have common parameters defining themselves simultaneously.}

%

\begin{definition}\label{def:com-paramter-condi-num}
Let $\bsA \bsX=\bsB$, where $\bsA$ $\in \R^{n \times n}$ \red{is an invertible matrix} whose entries are differentiable functions of a real parameter vector $\Omega_{\bsA }=[\bs\omega_{1}, \bs\omega_{2}, \ldots, \bs\omega_{p}, \bs\omega_{p+1} ,\ldots, \bs\omega_{N}] ^{\top} \in \R^{N}$, $\bsB$ $\in \R^{n \times m}$ whose entries are differentiable functions of a real parameter vector $\Omega_{\bsB}=[\bs\omega_{1}, \bs\omega_{2}, \ldots, \bs\omega_{p}, \bs\omega^\bsB_{p+1},\ldots,\bs\omega^\bsB_{M}] ^{\top}\in \R^{M}$, and $\mathbf{0}\neq \bsX \in \R^{n \times m}$. Let \red{the weight vectors corresponding to the parameterize representation  $\Omega_{\bsB }$ and $\Omega_{\bsA   }$ are given by}   $\mathbf{0}\leq \bsf  =[\bse_{1},\bse_{2}, \ldots, $ $\bse_{p}, \bsf_{p+1},\ldots, \bsf_{M} ]^\top \in \R^{M}$ and $\red{\mathbf{0}\leq \bsE}=[\bse_{1},\bse_{2}, \ldots, \bse_{p}, \bse_{p+1}, \ldots, \bse_{N} ]^\top \in \R^{N}$, \red{respectively}. Then, we define
\begin{align}\label{def:cond_parame}
\cond_{\bsE,\bsf }(\bsA(\Omega_\bsA),\bsB(\Omega_\bsB)):=
\lim_{\eta \rightarrow 0}\sup\Bigg  \{ \frac{1}{\eta}\frac{\|\delta \bsX\|_{\max}}{\|\bsX\|_{\max}}:\bsA(\Omega_\bsA+&\delta\Omega_\bsA)(\bsX+\delta\bsX)=\bsB(\Omega_{\bsB}+\delta \Omega_{\bsB} ),\nonumber\\
&|\delta \Omega_\bsA |\leq \eta\bsE, |\delta \Omega_{\bsB} |\leq \eta\bsf
\Bigg \}.
\end{align}
\end{definition}

In order to find an explicit formula for $\cond_{\bsE,\bsf }(\bsA(\Omega_\bsA),\bsB(\Omega_\bsB))$, we need the next Lemma \ref{lemm:parti-deri}, which states the explicit derivative expressions of $\bsA^{-1}$ and $\bsX$ with respect to the parameter vectors $\Omega_{\bsA}$ and $\Omega_{\bsB }$. Since the proof of the following lemma is trivial and we omit it here.

\begin{lemma}\label{lemm:parti-deri}
Let $\bsA \bsX=\bsB$, where $\bsA$ $\in \R^{n \times n}$ is an invertible matrix whose entries are differentiable functions of a vector of real parameters
$$
\Omega_{\bsA }=[ \bs\omega_{1}, \bs\omega_{2}, \ldots, \bs\omega_{p}, \bs\omega_{p+1} ,\ldots, \bs\omega_{N}] ^{\top} \in \R^{N},
$$
 \red{$\bsB \in \R^{n \times m}$ whose entries} are differentiable function of a vector of parameters
$$
\Omega_{\bsB  }=[ \bs\omega_{1}, \bs\omega_{2}, \ldots, \bs\omega_{p}, \bs\omega^\bsB_{p+1} ,\ldots, \bs\omega^\bsB_{M}] ^{\top} \in \R^{M},
$$
and $0\neq \bsX$ $\in \R^{n\times m}$. Then, the following equalities hold:
\begin{align*}
&{\rm (1)} \quad \frac{\partial \bsA^{-1}}{\partial \bs\omega_{k}}=-\bsA^{-1}\frac{\partial \bsA}{\partial\bs\omega_{k}}\bsA^{-1}, \ \quad\quad\quad\quad k=1,\ldots, p,\\
&{\rm (2)} \quad \frac{\partial \bsX}{\partial {\bs  \omega}^\bsB_{k}}=\bsA^{-1}   \frac{\partial \bsB }{\partial {\bs  \omega}^\bsB_{k}  } , \quad\quad\quad\quad\quad\quad\quad k=p+1,\ldots,M,\\
&{\rm (3)} \quad \frac{\partial\bsX}{\partial \bs\omega_{k}}=-\bsA^{-1}\frac{\partial \bsA}{\partial\bs\omega_{k}}\bsX,
\quad \quad\quad\quad\quad \quad k=p+1,\ldots,N, \\
&{\rm (4)} \quad \frac{\partial\bsX}{\partial \bs\omega_{k}}=-\bsA^{-1}\frac{\partial \bsA}{\partial\bs\omega_{k}}\bsX+\bsA^{-1}   \frac{\partial \bsB }{\partial {\bs  \omega}_{k}},  \quad k=1,\ldots,p.
\end{align*}
\end{lemma}

%



In the next theorem, we give the explicit expressions of the componentwise relative condition number $\cond_{\bsE,\bsf }(\bsA(\Omega_\bsA),\bsB(\Omega_\bsB))$ introduced in Definition \ref{def:com-paramter-condi-num}.
%

\begin{theorem}\label{th:explict-express-com-param-con-num1}
Let $\bsA \bsX=\bsB$, where $\bsA$ $\in \R^{n \times n}$ is an invertible matrix whose entries are differentiable functions of \red{ real parameter vector} \red{$\Omega_{\bsA }=[ \bs\omega_{1}, \ldots, \bs\omega_{p}, \bs\omega_{p+1} ,\ldots, \bs\omega_{N}] ^{\top} \in \R^{N}$}, $\bsB \in \R^{n \times m}$ whose entries are differentiable functions of \red{real parameter vector} \red{$\Omega_{\bsB  }=[ \bs\omega_{1}, \bs\omega_{2}, \ldots, \bs\omega_{p}, \bs\omega^\bsB_{p+1} ,\ldots, \bs\omega^\bsB_{M}] ^{\top} \in \R^{M}$}, and $0\neq \bsX$ $\in \R^{\red{n\times m}}$. \red{Let $\mathbf{0}\leq \bsf  =[\bse_{1},\bse_{2}, \ldots, $ $\bse_{p}, \bsf_{p+1},\ldots, \bsf_{M} ]^\top \in \R^{M}$ and $\bsE=[\bse_{1},\bse_{2}, \ldots, \bse_{p}, \bse_{p+1}, \ldots, \bse_{N} ]^\top \in \R^{N}$ with nonnegative entries}. Then \red{the explicit expression for $\cond_{\bsE,\bsf }(\bsA(\Omega_{A}),\bsB(\Omega_\bsB))$ is given by}
\begin{equation}\label{eq:Th c}
\cond_{\bsE,\bsf }(\bsA(\Omega_\bsA),\bsB(\Omega_\bsB))=c/ \|\bsX\|_{\max},
\end{equation}
where
$$
c=\Bigg \|\sum^{p}_{k=1}\left|\bsA^{-1}\frac{\partial \bsA}{\partial\bs\omega_{k}}\bsX-\bsA^{-1}   \frac{\partial \bsB }{\partial {\bs  \omega}_{k}  }\right | \bse_{k}+ \sum^{N}_{k=p+1}\left | \bsA^{-1}\frac{\partial \bsA}{\partial{ \bs \omega} _{k}}\bsX \right | \bse_{k}  +
\sum^{M}_{k=p+1} \left | \bsA^{-1}   \frac{\partial \bsB }{\partial {\bs  \omega}^\bsB_{k}  } \right| \bsf_{k} \Bigg \|_{\max}.
$$
\end{theorem}
\proof Because $\bsA$ and $\bsB$ are differentiable functions of the parameter \red{vectors} $\Omega_\bsA $ and  $\Omega_\bsB$, from $\bsX=\bsA^{-1}\bsB$, it is easy to see that $\bsX$ is also a differentiable function of $\Omega_\bsA $ and $\Omega_\bsB$, we can use differential calculus to obtain the following result:
\begin{align*}
\delta\bsX=\sum^{p}_{k=1}\frac{\partial\bsX}{\partial\bs\omega_{k}}\delta\bs\omega_{k}+ \sum^{N}_{k=p+1}\frac{\partial\bsX}{\partial\bs\omega_{k}}\delta\bs\omega_{k}+ \sum^{M}_{k=p+1}\frac{\partial\bsX}{\partial\bs\omega^\bsB_{k}}\delta\bs\omega^\bsB_{k} +\Oh(\|(\delta\Omega_\bsA,\delta \Omega_\bsB)\|^{2} ).
\end{align*}
From Lemma \ref{lemm:parti-deri},  it yields that
\begin{align*}
\delta\bsX&=\sum^{p}_{k=1}\left(-\bsA^{-1}\frac{\partial \bsA}{\partial\bs\omega_{k}}\bsX+\bsA^{-1}   \frac{\partial \bsB }{\partial {\bs  \omega}_{k}  }\right )\delta\bs\omega_{k}+ \sum^{N}_{k=p+1}\left (-\bsA^{-1}\frac{\partial \bsA}{\partial\bs\omega_{k}}\bsX \right )\delta\bs\omega_{k}\notag \\
&\quad\quad+
\sum^{M}_{k=p+1} \left (\bsA^{-1}   \frac{\partial \bsB }{\partial {\bs  \omega}^\bsB_{k}  } \right)\delta{\bs  \omega}^\bsB_{k} +\Oh(\|(\delta\Omega_\bsA,\delta\Omega_\bsB)\|^{2} ).
\end{align*}
Thus, \red{by} taking $\vect$ operation,  we obtain that
\begin{align}\label{eq:1}
\vect(\delta\bsX)&=\sum^{p}_{k=1}\vect \left(-\bsA^{-1}\frac{\partial \bsA}{\partial\bs\omega_{k}}\bsX+\bsA^{-1}   \frac{\partial \bsB }{\partial {\bs  \omega}_{k}  }\right )\delta\bs\omega_{k}+ \sum^{N}_{k=p+1}\vect \left(-\bsA^{-1}\frac{\partial \bsA}{\partial\bs\omega_{k}}\bsX\right )\delta\bs\omega_{k}\notag \\
&\quad \quad +
\sum^{M}_{k=p+1} \vect \left (\bsA^{-1}   \frac{\partial \bsB }{\partial {\bs  \omega}^\bsB_{k}  } \right )\delta{\bs  \omega}^\bsB_{k} +\Oh(\|(\delta\Omega_\bsA,\delta\Omega_\bsB)\|^{2}) \notag \\
&=\bsC \Diag (\bsa ) \frac{\delta \bsa}{\bsa }+\Oh(\|(\delta\Omega_\bsA,\delta\Omega_\bsB)\|^{2}), \,
\end{align}
where $\bsC \in \R^{(mn) \times (N+M-p)}$,  $\bsa\in \R^{N+M-1}$ and $\delta\bsa\in \R^{N+M-1}$ \red{are given by}
\begin{align*}
&	\bsC(:,1:p)=\left[\vect \left(-\bsA^{-1}\frac{\partial \bsA}{\partial\bs\omega_{1}}\bsX+\bsA^{-1}   \frac{\partial \bsB }{\partial {\bs  \omega}_{1}  }\right ), \ldots, \vect \left(-\bsA^{-1}\frac{\partial \bsA}{\partial\bs\omega_{p}}\bsX+\bsA^{-1}   \frac{\partial \bsB }{\partial {\bs  \omega}_{p}  }\right ) \right],\\
&	\bsC(:,(p+1):N)=\left[\vect \left(-\bsA^{-1}\frac{\partial \bsA}{\partial\bs\omega_{p+1}}\bsX\right ),\ldots, \vect \left(-\bsA^{-1}\frac{\partial \bsA}{\partial\bs\omega_{N}}\bsX\right )  \right],\\
&	\bsC(:,(N+1):( N+M-p))=\left[\vect \left (\bsA^{-1}   \frac{\partial \bsB }{\partial {\bs  \omega}^\bsB_{p+1}  } \right ),\ldots, \vect \left (\bsA^{-1}   \frac{\partial \bsB }{\partial {\bs  \omega}^\bsB_{M}  } \right ) \right ],\\
&\delta \bsa=[\delta\Omega_{\bsA}^\top ,\delta  \bs\omega^\bsB_{p+1} ,\ldots, \delta \bs\omega^\bsB_{M}]^\top,\quad \bsa=[\bsE^\top, \bsf_{p+1}, \ldots, \bsf_{M} ]^\top.
	\end{align*}
In the last equality\red{,} we use the fact that \red{if} $\bsa_i=0$ holds, then $\delta \bsa_i$ must be zero. Moreover,  the inequalities $|\delta \Omega_\bsA |\leq \eta\,\bsE$ and $ |\delta \Omega_{\bsB} |\leq \eta\bsf $ \red{imply}
$$
\left\| \frac{\delta \bsa}{\bsa } \right \|_\infty \leq \eta.
$$
From (\ref{eq:1}) and the  property  of infinity norm,   we deduce that
\begin{align*}
\left\|\delta\bsX\right\|_{\max}&=\left\|\vect(\delta\bsX)\right\|_{\infty} \leq \eta \left \|\bsC \Diag(\bsa) \right \|_\infty+\Oh(\|(\delta\Omega,\delta\bsB)\|^{2}),  \\
&= \eta \Big \| |\bsC| ~| \Diag(\bsa)|~ {\bf 1} \Big \|_\infty+\Oh(\big\|(\delta\Omega,\delta\bsB\big)\|^{2}),\\
&=\eta \Bigg \|\sum^{p}_{k=1}\left|\bsA^{-1}\frac{\partial \bsA}{\partial\bs\omega_{k}}\bsX-\bsA^{-1}   \frac{\partial \bsB }{\partial {\bs  \omega}_{k}  }\right | \bse_{k}+ \sum^{N}_{k=p+1}\left | \bsA^{-1}\frac{\partial \bsA}{\partial { \bs \omega}_k }\bsX \right | \bs\omega_{k},  \\
&\quad \quad +
\sum^{M}_{k=p+1} \left | \bsA^{-1}   \frac{\partial \bsB }{\partial {\bs  \omega}^\bsB_{k}  } \right| \bsf_{k} \Bigg \|_{\max}+\Oh(\|(\delta\Omega_\bsA,\delta\Omega_\bsB)\|^{2}),
\end{align*}
where $\bf 1$ is a vector with all components being one.
Let $\eta$ tend to zero, from Definition \ref{def:com-paramter-condi-num},  then \red{we finish the proof of this theorem}.
\eproof


When the entries of multiple right-hand sides  $\bsB$ are independent of the parameter vector $\Omega_{\bsA}$ of the coefficient matrix $\bsA$, we consider the case that $\bsB$ is a function of its entries $\bsb_{i,j}$. Suppose that $|\delta \bsB| \leq \eta \,\bsF$ in Definition \ref{def:com-paramter-condi-num}. It is easy to see that
\begin{equation}\label{eq:unB}
	\frac{\partial \bsB}{\partial\bsb_{i,j}}=e_{i,(n)}e_{j,(m)}^{\top}, \quad \quad i=1,\ldots,n, \quad j=1,\ldots,m,
\end{equation}
where $e_{i,(n)} \in \R^{n}$ is the $i$-th column of the $n$-by-$n$ identity matrix $\bsI_n$. Since there \red{are} no common parameters of $\Omega_{\bsA}$ and $B$, we have $p=0$ in the expression of $\cond_{\bsE,\bsF}(\bsA(\Omega_{\red{A}}),\bsB)$ given by Theorem \ref{th:explict-express-com-param-con-num1}. Substituting \eqref{eq:unB} into \eqref{eq:Th c} and noting
$$
\sum_{i=1}^n \sum_{j=1}^m \left | \bsA^{-1} e_{i,(n)}e_{j,(m)}^{\top} \right | |\bsf_{i,j}|=\left| \bsA^{-1}\right| |\bsF|,
$$
we have the following theorem.
%
%
%


\begin{theorem}\label{th:explict-express-com-param-con-num2}
Let $\bsA \bsX=\bsB$, where $\bsA$ $\in \R^{n \times n}$ is an invertible matrix whose entries are differentiable functions of a vector of real parameters \red{$\Omega_{\bsA }=[ \bs\omega_{1}, \bs\omega_{2}, \ldots, \bs\omega_{p}, \bs\omega_{p+1} ,\ldots, \bs\omega_{N}] ^{\top} \in \R^{N}$}, $\bsB \in \R^{n \times m}$ whose parameter \red{vector} is the set of its all entries, and $0\neq \bsX$ $\in \R^{\red{n \times m}}$. \red{Let $\mathbf{0}\leq\bsE=[\bse_{1},\bse_{2}, \ldots, \bse_{p}, \bse_{p+1}, \ldots, \bse_{N} ]^\top \in \R^{N}$ and $\mathbf{0}\leq \bsF \in \R^{n\times m}$. In addition, we assume that $\Omega_\bsA \cap \bsB= \emptyset$}. For $\cond_{\bsE,\bsF}(\bsA(\Omega_{\bsA }),\bsB )$ defined in Definition \ref{def:com-paramter-condi-num}, we have its explicit expression as \red{follows}
$$
\cond_{\bsE,\bsF}(\bsA(\Omega_\bsA),\bsB)=\frac{\Big\||\bsA^{-1}|\bsF+\sum^{N}\limits_{k=1}\big|\bsA^{-1}\frac{\partial\bsA}{\partial\bs\omega_{k}}\bsX\big|\bse_{k}\Big\|_{\max}}{\|\bsX\|_{\max}}.
$$
\end{theorem}

In the following, we will show the explicit expression of $\cond_{\bsE,\bsF}(\bsA,\bsB )$ in \eqref{th:the expression-compon} can be deduced from Theorem \ref{th:explict-express-com-param-con-num2} by considering $\Omega_\bsA $ as the set of entries of $\bsA$. Since
\begin{align*}
\cond_{\bsE,\bsF}(\bsA,\bsB )=\frac{\Big\|\big|\bsA^{-1}\big|\bsF+ \sum^{n}\limits_{i=1}\sum^{m}\limits_{j=1}\left|\bsA^{-1}\frac{\partial\bsA}{\partial\bsa_{i,j}}\bsX\right|\bse_{i,j}\Big\|_{\max}}{\|\bsX \|_{\max}},
\end{align*}
 and using the fact that $\frac{\partial\bsA}{\partial \,\bsa_{i,j}}=e_{i}e_{j}^{\top}$, we have
\begin{align*}
\sum^{n}\limits_{i=1}\sum^{n}\limits_{j=1}\big|\bsA^{-1}\frac{\partial\bsA}{\partial\bsa_{i,j}}\bsX\big|\bse_{i,j}&=\sum^{n}_{i=1}\sum^{n}_{j=1}\Big|\bsA^{-1}e_{i}e_{j}^{\top}\bsX\Big|\bse_{i,j}=\sum^{n}_{i=1}\sum^{n}_{j=1}\Big|\bsA^{-1}(:,i)\bsX(j,:)\Big|\bse_{i,j},\\
&=\sum^{n}_{j=1}\left|\sum^{n}_{i=1}\bsA^{-1}(:,i)\bse_{i,j}\right|\big|\bsX(j,:)\big| =|\bsA^{-1}|\bsE|\bsX|,
\end{align*}
which gives the expression of $\cond_{\bsE,\bsF}(\bsA,\bsB)$ in \eqref{th:the expression-compon}.

 \red{As shown in \cite{Mary2017,Amestoy2015f,Amestoy2017f}, the sparsity of right-hand sides has practical background and enhance the efficiency of the direct solver for \eqref{eq:mul}. Therefore, we introduce the sparse representation  of $\bsB$  in \eqref{eq:b}, where we assume that  the non-zero entries of $\bsB$  are independent of the parameters $\Omega_{\bsA }$ defining $\bsA$. Denote}
\begin{equation}\label{eq:b}
	\bsB=\sum_{k=1}^M  {\bs  \omega}^\bsB_{k} \bsS_k,
\end{equation}
where $ \bsS_k \in \R^{n\times m}$ are matrices of constants, typically 0s and 1s, which \red{describe} the sparse pattern of $\bsB$ and \red{$\Omega_{B}=\{{\bs\omega}^\bsB_{k}\}^{M}_{k=1}$} are independent parameters. In the reminder of this paper, we always assume that the parameters \red{$\Omega_{B}=\{{\bs\omega}^\bsB_{k}\}^{M}_{k=1}$} are independent of \red{the parameter vector} $\Omega_{\bsA }$ defining the matrix $\bsA$. Thus from Lemma \ref{lemm:parti-deri} and \eqref{eq:b}, we deduce that
$$
\frac{\partial \bsB}{\partial {\bs  \omega}^\bsB_{k}}= \bsS_k.
$$

Noting the above equation, the following theorem holds from Theorem \ref{th:explict-express-com-param-con-num1}.

\begin{theorem}\label{th:explict-express-com-param-con-num}
	Let $\bsA \bsX=\bsB$. \red{Assume that $\bsB$ has the sparse representation \eqref{eq:b} and}	 $\bsA$ $\in \R^{n \times n}$ is an invertible matrix whose entries are differentiable functions of real parameter vector \red{$\Omega_{\bsA }=[\bs\omega_{1}, \bs\omega_{2}, \ldots, \bs\omega_{p}, \bs\omega_{p+1} ,\ldots, \bs\omega_{N}] ^{\top} \in \R^{N}$}, $\bsB \in \R^{n \times m}$ whose entries are differentiable functions of \red{ $\Omega_{B}=\{{\bs\omega}^\bsB_{k}\}^{M}_{k=1}$. Suppose that} $0\neq \bsX$ $\in \R^{n\times m}$. \red{Let $\mathbf{0}\leq\bsE=[\bse_{1},\bse_{2}, \ldots, \bse_{N} ]^\top \in \R^{N}$ and $\mathbf{0}\leq \bsf \in \R^{M}$, and assume that $\Omega_\bsA \cap\,\Omega_\bsB= \emptyset$ and $\bsE \cap\,\bsf= \emptyset$}. For $\cond_{\bsE,\bsf}(\bsA(\Omega_\bsA),\bsB(\Omega_\bsB))$ defined in Definition \ref{def:com-paramter-condi-num}, we have its explicit expression as \red{follows}
$$
\cond_{\bsE,\bsf }(\bsA(\Omega_\bsA), \bsB(\Omega_\bsB))=\frac{\Bigg \| \sum^{N}\limits_{k=1}\left | \bsA^{-1}\frac{\partial \bsA}{\partial{ \bs \omega} _{k}}\bsX \right | \bse_{k}  +
\sum^{M}\limits_{k=1} \left | \bsA^{-1}  \bsS_k \right| \bsf_{k} \Bigg \|_{\max}} {\|\bsX\|_{\max}  }.
$$
\end{theorem}

\begin{remark}
 If $\bsB$ is sparse, its sparse patten is persevered during numerical computations since only non-zero entries and its positions are stored in computer.Thus it is reasonable to measure the conditioning of the linear system with multiple right-hand sides \eqref{eq:mul} under componentwise perturbation analysis. If $\bsB$ is not sparse \red{and hence is} general unstructured, i.e., the parameter representation of $\bsB$ is \red{just} the set of all entries of $\bsB$, we can take $\bsS_{ij}=e_{i,(n)}e_{j,(m)}^{\top}$, thus Theorem \ref{th:explict-express-com-param-con-num} reduces to Theorem \ref{th:explict-express-com-param-con-num2}.
\end{remark}

In the next corollary, we will consider componentwise condition numbers of \eqref{eq:mul} for the case that $\bsA$ is a general unstructured matrix and $\bsB$ is sparse.\,\,In this situation, the perturbation on $\bsA$ is assumed to be arbitrary while the perturbation on $\bsB $ should preserve the sparse pattern of $\bsB$.

\begin{definition}\label{def:com-paramter-condi-num1}
Let $\bsA \bsX=\bsB$, where $\bsA$ $\in \R^{n \times n}$ is a nonsingular matrix whose parameter representation vector $\Omega_{\bsA}$ is the set of all entries of $\bsA$, $\bsB$ $\in \R^{n \times m}$ whose entries are differentiable function of a vector of parameters $\Omega_{\bsB }=[ \bs\omega_{1}^\bsB ,\ldots, \bs\omega^\bsB_{M}] ^{\top}$ $\in \R^{M}$ such that $\bsB$ can be expressed as \eqref{eq:b}, and $\mathbf{0}\neq \bsX$ $\in \R^{n \times m}$. Let $\mathbf{0}\leq \bsf  =[\bsf_{1}, \ldots, \bsf_{M} ]^\top \in \R^{M}$ and $\bsE \in \R^{n\times n}$ with nonnegative entries. Then, we define
\begin{align*}
\cond_{\bsE,\bsf }(\bsA,\bsB(\Omega_\bsB)):=
\lim_{\eta \rightarrow 0}\sup\Bigg  \{ \frac{1}{\eta}\frac{\|\delta \bsX\|_{\max}}{\|\bsX\|_{\max}}: (\bsA+\delta \bsA)(\bsX+\delta\bsX)=\bsB(\Omega_{\bsB}+\delta \Omega_{\bsB} ),  \\
|\delta \bsA  |\leq \eta\,\bsE, |\delta \Omega_{\bsB} |\leq \eta\bsf
\Bigg \}.
\end{align*}
\end{definition}

\begin{corollary}\label{co:3.1}
	Let $\bsA \bsX=\bsB$, where $\bsA$ $\in \R^{n \times n}$ is an invertible matrix and $\bsB \in \R^{n \times m}$ \red{is a sparse matrix which has the sparse representation  \eqref{eq:b}}, and $0\neq \bsX$ $\in \R^{n \times m}$. \red{Let $\mathbf{0}\leq \bsf  =[\bsf_{1}, \ldots, \bsf_{M} ]^\top \in \R^{M}$ and $\bsE \in \R^{n\times n}$ with nonnegative entries}. For $\cond_{\bsE,\bsf}(\bsA,\bsB(\Omega_\bsB))$ defined in Definition \red{\ref{def:com-paramter-condi-num1}}, we have its explicit expression as \red{follows}:
\begin{equation*}
\cond_{\bsE, \bsf }(\bsA,\bsB(\Omega_\bsB))=\frac{\Bigg \| \left|\bsA^{-1}\right |\, \bsE\, |\red{\bsX}|  +
\sum^{M}\limits_{k=1} \left | \bsA^{-1}  \bsS_k \right| \bsf_{k} \Bigg \|_{\max}} {\|\bsX\|_{\max}  }.
\end{equation*}
\red{
Furthermore, if $\bsE=|A|$ and $\bsf=|\Omega_\bsB| $ is given in \eqref{eq:b},  we have
\begin{equation}\label{eq:cond_unstruWithSparseRH}
\cond_{|\bsA|, |\Omega_\bsB| }(\bsA,\bsB(\Omega_\bsB))=\frac{\Bigg \| \left|\bsA^{-1}\right |\, |\bsA|\, |\red{\bsX}|  +
\sum^{M}\limits_{k=1} \left | \bsA^{-1}  \bsS_k \right| {|\bs \omega_{k}^{\bsB }|} \Bigg \|_{\max}} {\|\bsX\|_{\max}  }.
\end{equation}
}
\end{corollary}

In the following section, we will apply Theorem \ref{th:explict-express-com-param-con-num} to the case that the coefficiment matrix $\bsA$ in \eqref{eq:mul} belongs to \{1;1\}-quasiseparable matrices and multiple right-hand sides  $B$ is sparse.

\section{Condition number for  multiple right-hand side linear systems with  \{1;1\}-quasiseparable coefficient matrices }\label{section:quasi-represetation}

In this section, when the matrix $\bsB$ of \red{the multiple right-hand side linear system} \eqref{eq:mul} has the sparse representation  \eqref{eq:b} and $\bsA$ is a  \{1;1\}-quasiseparable matrix, we will focus on deriving explicit expressions for the condition number of the multiple right-hand side  linear system with respect to the quasiseparable representation and the Givens-vector representation via tangent.

\subsection{The general quasiseparable representation}
In this subsection, we will consider the componentwise condition number for the solution of multiple right-hand side linear system \eqref{eq:mul}, \red{where $\bsB$ has the sparse representation \eqref{eq:b}}, and $\bsA$ has the quasiseparable representation $\Omega_{QS}$ given in Definition \ref{th:{1;1}}. Since the norm of the vector $\Omega_{QS}$ of parameters does not determine the norm of the matrix $\bsA$, it is more natural to consider componentwise perturbations of $\Omega_{QS}$. We will apply Theorem \ref{th:explict-express-com-param-con-num} to derive the explicit expressions for the componentwise condition number of the solution of \red{the multiple right-hand side linear system} \eqref{eq:mul} with respect to the general quasiseparable representation $\Omega_{QS}$ \red{given in Definition \ref{th:{1;1}}}.


 \red{Recall that when $\bsA$ has a general quasiseparable representation \eqref{eq:qs general} and $\bsB \in \R^{n\times m}$ \red{has the sparse representation} \eqref{eq:b}, in view of Definition \ref{def:com-paramter-condi-num}, we introduce  the weight vector
 \begin{equation}\label{eq:eps}
\bsE_{QS}=\Big(\{\bse_{\bsp_{i}}\}^{n}_{i=2}, \{\bse_{\bsa_{i}}\}^{n-1}_{i=2}, \{\bse_{\bsq_{i}}\}^{n-1}_{i=1}, \{\bse_{\bsg_{i}}\}^{n-1}_{i=1}, \{\bse_{\bsb_{i}}\}^{n-1}_{i=2}, \{\bse_{\bsh_{i}}\}^{n}_{i=2}, \{\bse_{\bsd_{i}}\}^{n}_{i=1}\Big)
 \end{equation}
 with respect to the general quasiseparable representation \eqref{eq:qs general}  of $\bsA$. Therefore, using Definition \ref{def:com-paramter-condi-num}, we introduce the structure componentwise condition number for \eqref{eq:mul} with respect to the general quasiseparable representation \eqref{eq:qs general}  of $\bsA$ and the sparse representation \eqref{eq:b} of $\bsB$ in Theorem \ref{th:important- expres-compon}.
}
\begin{theorem}\label{th:important- expres-compon}
Let $\bsA \bsX=\bsB$, where $\mathbf{0}\neq \bsX \in \R^{n \times m}$, $\bsA$ $\in \R^{n \times n}$ is a nonsingular \{1;1\}-quasiseparable matrix with a quasiseparable representation $\Omega_{QS}$ \eqref{eq:qs general}, $\bsB \in \R^{n\times m}$ \red{has the sparse representation} \eqref{eq:b}, and $\bsA=\bsA_{\bsL}+\bsA_{\bsD}+\bsA_{\bsU}$, with $\bsA_{\bsL}$ strictly lower triangular, $\bsA_{\bsD}$ diagonal, and $\bsA_{\bsU}$ strictly upper triangular. Let the weight vectors  $\mathbf{0}\leq\bsf \in \R^{M} $ and $\mathbf{0}\leq \bsE_{QS}$ $\in \R^{7n-8}$, where $\bsE_{QS}$ is given by \eqref{eq:eps}. Suppose $\Omega_{QS} \cap \Omega_\bsB= \emptyset$. \red{Then}
\begin{align*}
\cond_{\bsE_{QS},\bsf }(\bsA(\Omega_{QS}),\bsB(\Omega_\bsB))&=\frac{1}{\|\bsX\|_{\max}}\Bigg\| \sum^{M}_{k=1} \left | \bsA^{-1}  \bsS_k \right| \bsf_{k}   +  |\bsA^{-1}||\bsQ_{\bsd}||\bsX| + |\bsA^{-1}||\bsQ_{\bsp}||\bsA_{\bsL}\bsX|\\
&\quad +|\bsA^{-1}\bsA_{\bsL}||\bsQ_{\bsq}||\bsX| + |\bsA^{-1}\bsA_{\bsU}||\bsQ_{\bsg}||\bsA_{\bsU}\bsX| +  |\bsA^{-1}\bsA_{\bsU}||\bsQ_{\bsh}||\bsX| \\
&\quad+\sum^{n-1}_{i=2}\left| \bsA^{-1}\begin{bmatrix}0&0\\\bsA(i+1:n,1:i-1)&0\\\end{bmatrix}\bsX\right| \left|\frac{\bse_{\bsa_{i}}}{\bsa_{i}}\right| \\
&\quad+\sum^{n-1}_{i=\red{2}}\left| \bsA^{-1} \begin{bmatrix} 0&\bsA(1:i-1,i+1:n)\\0&0\\ \end{bmatrix} \bsX\right| \left|\frac{\bse_{\bsb_{i}}}{\bsb_{i}}\right|  \Bigg\|_{\max},
\end{align*}
where
\begin{align} \label{eq:qsE}
    \Omega_{QS}&=\Big(\{\bsp_{i}\}^{n}_{i=2}, \{\bsa_{i}\}^{n-1}_{i=2}, \{\bsq_{i}\}^{n-1}_{i=1}, \{\bsg_{i}\}^{n-1}_{i=1}, \{\bsb_{i}\}^{n-1}_{i=2}, \{\bsh_{i}\}^{n}_{i=2}, \{\bsd_{i}\}^{n}_{i=1}\Big),\nonumber \\
    \bsQ_{\bsd}&=\diag\Big(\bse_{\bsd_{1}}, \ldots, \bse_{\bsd_{n}}\Big),   \bsQ_{\bsp}=\diag\Big(1, \frac{\bse_{\bsp_{2}}}{\bsp_{2}}, \ldots, \frac{\bse_{\bsd_{n}}}{\bsp_{n}}\Big),\notag\\
    \bsQ_{\bsq}&=\diag\Big(\frac{\bse_{\bsq_{1}}}{\bsq_{1}}, \ldots, \frac{\bse_{\bsq_{n-1}}}{\bsq_{n-1}}, 1\Big),
    \bsQ_{\bsg}=\diag\Big(\frac{\bse_{\bsg_{1}}}{\bsg_{1}}, \ldots, \frac{\bse_{\bsg_{n-1}}}{\bsg_{n-1}}, 1\Big),\notag\\
    \bsQ_{\bsh}&=\diag\Big(1, \frac{\bse_{\bsh_{2}}}{\bsh_{2}}, \ldots, \frac{\bse_{\bsh_{n}}}{\bsh_{n}}\Big).  
\end{align}
\end{theorem}
\proof The proof is proceed by calculating the contribution of each subset of parameters to the expression for $\cond_{\bsE_{QS},\bsf }(\bsA(\Omega_{QS}),\bsB(\Omega_\bsB))$ given in Theorem \ref{th:explict-express-com-param-con-num} step by step.  For the parameters $\{\bsd_{i}\}^{n}_{i=1}$, from 1) in Lemma \ref{lemm:proper-all-param} together with Theorem \ref{th:explict-express-com-param-con-num}, using the fact
\begin{equation}\label{eq:out}
	\left|\bsa \cdot  \bsb^\top\right|=\big|\bsa \big| \cdot \left| \bsb^\top \right|,
\end{equation}
where $\bsa$ and $\bsb$ are two column vectors, we derive that
\begin{align*}
\kappa_{\bsd}:&= \sum^{n}_{i=1} \Big|\bsA^{-1}\frac{\partial\bsA}{\partial\bsd_{i}}\bsX\Big||\bse_{\bsd_{i}}|=\sum^{n}_{i=1}\Big|\bsA^{-1}e_{i}e^{\top}_{i}\bsX \Big||\bse_{\bsd_{i}}| =\sum^{n}_{i=1}\big|\bsA^{-1}e_{i}\big|\big|e^{\top}_{i}\bsX \big||\bse_{\bsd_{i}}|=|\bsA^{-1}||\bsQ_{\bsd}||\bsX|.
\end{align*}
 As for the parameters $\{\bsp_{i}\}^{n}_{i=2}$, from 2) in Lemma \ref{lemm:proper-all-param} and using \eqref{eq:out} again, it can be deduced that
\begin{align*}
\kappa_{\bsp}:&=\sum^{n}_{i=1} \Big|\bsA^{-1}\frac{\partial\bsA}{\partial\bsp_{i}}\bsX\Big|\big|\bse_{\bsp_{i}}\big|=\sum^{n}_{i=1}\left|\bsA^{-1}\bsp_{i}\frac{\partial\bsA}{\partial\bsp_{i}}\bsX \right|\left|\frac{\bse_{\bsp_{i}}}{\bsp_{i}}\right|\\
&=\sum^{n}_{i=1}\left|\bsA^{-1}e_{i}\bsA_{\bsL}(i,:)\bsX \right|\left|\frac{\bse_{\bsp_{i}}}{\bsp_{i}}\right|
=\sum^{n}_{i=1}\left|\bsA^{-1}e_{i}||\bsA_{\bsL}(i,:)\bsX \right|\left|\frac{\bse_{\bsp_{i}}}{\bsp_{i}}\right|\\
&=\sum^{n}_{i=1}\left|\bsA^{-1}(:,i)||\bsA_{\bsL}\bsX(i,:) \right|\left|\frac{\bse_{\bsp_{i}}}{\bsp_{i}}\right|=|\bsA^{-1}||\bsQ_{\bsp}||\bsA_{\bsL}\bsX|.
\end{align*}
Now we want to computer the term related to  $\{\bsa_{i}\}^{n-1}_{i=2}$ in the condition number expression.  From 3) in Lemma \ref{lemm:proper-all-param}, we have
\begin{align*}
\kappa_{\bsa}:&=\sum^{n-1}_{i=2}\left|\bsA^{-1}\frac{\partial\bsA}{\partial\bsa_{i}}\bsX\right||\bse_{\bsa_{i}}|=\sum^{n-1}_{i=2}\left|\bsA^{-1}\bsa_{i}\frac{\partial\bsA}{\partial\bsa_{i}}\bsX\right|\left|\frac{\bse_{\bsa_{i}}}{\bsa_{i}}\right|\\
&=\sum^{n-1}_{i=2}\left| \bsA^{-1}\begin{bmatrix} 0&0\\\bsA(i+1:n,1:i-1)&0\\
\end{bmatrix} \bsX\right|\bigg|\frac{\bse_{\bsa_{i}}}{\bsa_{i}}\bigg|.
\end{align*}
Again for  $\{\bsq_{j}\}^{n-1}_{j=1}$, from $4)$ in Lemma \ref{lemm:proper-all-param} and \eqref{eq:out},  we have
\begin{align*}
\kappa_{\bsq}:&=\sum^{n-1}_{j=1} \Big|\bsA^{-1}\frac{\partial\bsA}{\partial\bsq_{j}}\bsX \Big||\bse_{\bsq_{j}}| = \sum^{n-1}_{j=1} \Big|\bsA^{-1}\bsA_{\bsL}(:,j)e^{\top}_{j}\bsX \Big| \left|\frac{\bse_{\bsq_{j}} }{\bsq_j }\right| \\
&= \sum^{n-1}_{j=1} \Big|\bsA^{-1}\bsA_{\bsL}(:,j)||\bsX(j,:) \Big||\bse_{\bsq_{j}}|
              =|\bsA^{-1}\bsA_{\bsL}||\bsQ_{\bsq}||\bsX|.
              \end{align*}
Similar for the contributing terms in  $\cond_{\bsE_{QS},\bsf }(\bsA(\Omega_{QS}), \bsX)$ with respect to parameters $\{\bsg_{i}\}^{n-1}_{\red{i=1}}$, $\{\bsb_{i}\}^{n-1}_{i=2}$, and $\{\bsh_{i}\}^{\red{n}}_{i=2}$, based on the above technique, we can derive that
\begin{align*}
\kappa_{\bsg}&:= \sum^{n-1}_{i=1}\Big| \bsA^{-1}\frac{\partial \bsA}{\partial\bsg_{i}}\bsX \Big|\Big|\frac{\bse_{\bsg_{i}}}{\bsg_{i}}\Big| = |\bsA^{-1}||\bsQ_{\bsg}||\bsA_{\bsU}\bsX|,\\
\kappa_{\bsh}&:=\sum^{n}_{i=2}\Big| \bsA^{-1}\frac{\partial\bsA}{\partial\bsh_{i}}\bsX\Big||\bse_{\bsh_{i}}| = |\bsA^{-1}\bsA_{\bsU}||\bsQ_{\bsh}||\bsX|,\\
\kappa_{\bsb}&:=\sum^{n-1}_{\red{i=2}}\left| \bsA^{-1}\frac{\partial\bsA}{\bsb_{i}}\bsX \right|=\sum^{n-1}_{\red{i=2}}\left| \bsA^{-1}\begin{bmatrix}  0&\bsA(1:i-1,i+1:n)\\0&0\\
\end{bmatrix} \bsX\right|\left|\frac{\bse_{\bsb_{i}}}{\bsb_{i}}\right|.
\end{align*}
According to Theorem \ref{th:explict-express-com-param-con-num} we have that
\begin{align*}
\cond_{\bsE_{QS},\bsf }(\bsA(\Omega_{QS}),\bsB(\Omega_{\bsB }))=\frac{1}{\|\bsX\|_{\max}}\Bigg\|        \sum^{M}_{k=1} \left | \bsA^{-1}  \bsS_k \right| \bsf_{k}   +\kappa_{\bsd}+\kappa_{\bsp}+\kappa_{\bsq}+\kappa_{\bsa}+\kappa_{\bsh}+\kappa_{\bsg}+\kappa_{\bsb} \Bigg\|_{\max},
\end{align*}
which completes the proof of this theorem. \eproof

%



\begin{remark}\label{remark:1}
If the multiple right-hand sides $\bsB$ is not sparse \red{but} general unstructured, from Theorem \ref{th:explict-express-com-param-con-num2}, the explicit expression for $\cond_{\bsE_{QS},\bsF}(\bsA(\Omega_{QS}),\bsB)$ of \eqref{eq:mul}, where $\bsA$ is a \{1;1\}-quasiseparable nonsingular matrix with a quasiseparable representation $\Omega_{QS}$ and $\bsB$ is a general dense matrix, can be characterized by
\begin{align}\label{eq:rm}
\cond_{\bsE_{QS},\bsF}(\bsA(\Omega_{QS}),\bsB )&\red{:}=\lim_{\eta \rightarrow 0}\sup\Bigg  \{ \frac{1}{\eta}\frac{\|\delta \bsX\|_{\max  }}{\|\bsX\|_{\max  }}:\bsA(\Omega_{QS} +\delta\Omega_{QS} )(\bsX+\delta\bsX)=\bsB+\delta\bsB,  \nonumber \\
&\hspace{4.2cm}  |\delta \Omega_{QS} |\leq \eta \,\bsE_{QS}, |\delta \bsB |\leq \eta\bsF
\Bigg \} \nonumber \\
&=\frac{1}{\|\bsX\|_{\max}}\Bigg\|  |\bsA^{-1}|\bsF  +  |\bsA^{-1}||\bsQ_{\bsd}||\bsX| + |\bsA^{-1}||\bsQ_{\bsp}||\bsA_{\bsL}\bsX| \nonumber \\
&\quad +|\bsA^{-1}\bsA_{\bsL}||\bsQ_{\bsq}||\bsX| + |\bsA^{-1}\bsA_{\bsU}||\bsQ_{\bsg}||\bsA_{\bsU}\bsX| +  |\bsA^{-1}\bsA_{\bsU}||\bsQ_{\bsh}||\bsX| \nonumber \\
&\quad+\sum^{n-1}_{i=2}\left| \bsA^{-1}\begin{bmatrix}0&0\\\bsA(i+1:n,1:i-1)&0\\\end{bmatrix}\bsX\right| \left|\frac{\bse_{\bsa_{i}}}{\bsa_{i}}\right| \nonumber \\
&\quad+\sum^{n-1}_{i=1}\left| \bsA^{-1} \begin{bmatrix} 0&\bsA(1:i-1,i+1:n)\\0&0\\ \end{bmatrix} \bsX\right| \left|\frac{\bse_{\bsb_{i}}}{\bsb_{i}}\right|  \Bigg\|_{\max}. 
\end{align}
For the linear system $\bsA \bsx =\bsb $ with a \{1;1\}-quasiseparable coefficient matrix $\bsA(\Omega_{QS})$ in the general quasiseparable representation $\Omega_{QS}$, Dopico and Pom{\'e}s  \cite{Dopico2016}  introduced the componentwise condition number of $\bsx$ with respect to  the general quasiseparable representation $\Omega_{QS}$  and the corresponding explicit formula  \cite[Theorem 4.5]{Dopico2016} is given below
\begin{align}\label{eq:rm1}
{\rm cond}_{\bsE_{QS},\bsy  }(\bsA(\Omega_{QS}),\bsb)&\red{:}=\lim_{\eta \rightarrow 0}\sup\Bigg  \{ \frac{1}{\eta}\frac{\|\delta \bsx\|_{\infty }}{\|\bsx\|_{\infty }}:\bsA(\Omega_{QS} +\delta\Omega_{QS} )(\bsx+\delta\bsx)=\bsb+\delta\bsb,  \nonumber \\
&\hspace{4.2cm}  |\delta \Omega_{QS} |\leq \eta \,\bsE_{QS}, |\delta \bsb|\leq \eta\bsy
\Bigg \} \nonumber \\
&=\frac{1}{\|\bsx\|_{\infty }}\bigg\|  |\bsA^{-1}|\bsy    +  |\bsA^{-1}||\bsQ_{\bsd}||\bsx| + |\bsA^{-1}||\bsQ_{\bsp}||\bsA_{\bsL}\bsx| \nonumber \\
&\quad+|\bsA^{-1}\bsA_{\bsL}||\bsQ_{\bsq}||\bsx| + |\bsA^{-1}\bsA_{\bsU}||\bsQ_{\bsg}||\bsA_{\bsU}\bsx| +  |\bsA^{-1}\bsA_{\bsU}||\bsQ_{\bsh}||\bsx| \nonumber \\
&\quad+\sum^{n-1}_{i=2}\left| \bsA^{-1}\begin{bmatrix}0&0\\\bsA(i+1:n,1:i-1)&0\\\end{bmatrix}\bsx\right| \Big|\frac{\bse_{\bsa_{i}}}{\bsa_{i}}\Big| \nonumber \\
&\quad+\sum^{n-1}_{i=1}\left| \bsA^{-1} \begin{bmatrix} 0&\bsA(1:i-1,i+1:n)\\0&0\\ \end{bmatrix} \bsx\right| \Big|\frac{\bse_{\bsb_{i}}}{\bsb_{i}}\Big|  \bigg\|_{\infty }, 
\end{align}
where $\Omega_{QS},\,  \bsE_{QS}, \, \bsQ_{\bsd},\,  \bsQ_{\bsq},\,   \bsQ_{\bsg}$ and $\bsQ_{\bsh}$ are defined in \eqref{eq:qsE}. When the column dimension of the multiple right-hand side linear system \eqref{eq:mul} is 1, \eqref{eq:mul} degenerates to the linear system $\bsA \bsx =\bsb$.\,\,Thus it is easy to see that \red{$\cond_{\bsE_{QS},\bsF}(\bsA(\Omega_{QS}),\bsB )$} given in \eqref{eq:rm} reduces to \red{${\rm cond}_{\bsE_{QS},\bsy  }(\bsA(\Omega_{QS}),\bsb)$} given by \eqref{eq:rm1} when $\bsF=\bsy$.

\red{In Theorem \ref{th:important- expres-compon}, it is natural to take the choice of $\bsE_{QS}$ \eqref{eq:eps} as $\bsE_{QS} =|\Omega_{QS}|$, where  $\Omega_{QS}$ is given by \eqref{eq:qs general} and is the general quasiseparable representation introduced in Definition \ref{th:{1;1}} of $\bsA$}. Under this situation, noting $\bsQ_{\bsd}=|\bsA_{\bsD}|$, $\bsQ_{\bsp}=\bsQ_{\bsq}=\bsQ_{\bsg}=\bsQ_{\bsh}=\bsI$, and $\big|\frac{\bse_{\bsa_{i}}}{\bsa_{i}}\big|=\big|\frac{\bse_{\bsb_{i}}}{\bsb_{i}}\big|=1$, we have
\red{
\begin{corollary}\label{th:important-express-cond_f}
	Let $\bsA\bsX=\bsB$, where $\mathbf{0}\neq\bsX \in\R^{n\times m}$, $\bsB \in \R^{n\times m}$ \red{has the sparse representation} \eqref{eq:b},  and $\bsA$ $\in \R^{n\times n}$ is a \{1;1\}-quasiseparable  nonsingular matrix \red{with a quasiseparable representation $\Omega_{QS}$ \eqref{eq:qs general}} such that $\bsA=\bsA_{\bsL}+\bsA_{\bsD}+\bsA_{\bsU}$, with $\bsA_{\bsL}$ strictly lower triangular, $\bsA_{\bsD}$ diagonal, $\bsA_{\bsU}$ strictly upper triangular. Then
\begin{align*}
\red{\cond_{|\Omega_{QS}|,|\Omega_\bsB|}(\bsA(\Omega_{QS}),\bsB(\Omega_\bsB) )} &= \frac{1}{\|\bsX\|_{\max}}\bigg\|\sum^{M}_{k=1} \left | \bsA^{-1}  \bsS_k \right| \red{|\bsw^{\bsB}_{k}|}  +|\bsA^{-1}||\bsA_{\bsD}||\bsX|
|\bsA^{-1}\bsA_{\bsL}||\bsX|\\
&\quad +|\bsA^{-1}||\bsA_{\bsL}\bsX|+|\bsA^{-1}||\bsA_{\bsU}\bsX|+|\bsA^{-1}\bsA_{\bsU}||\bsX|\\
&\quad +
\sum^{n-1}_{i=2}\Big|\bsA^{-1}
\begin{bmatrix}0&0\\\bsA(i+1:n,1:i-1)&0\\
\end{bmatrix}\bsX\Big|\\
&\quad+\sum^{n-1}_{i=2}\left|\bsA^{-1}
\begin{bmatrix}0&\bsA(1:i-1,i+1:n)\\0&0\\
\end{bmatrix}\bsX\right|  \bigg\|_{\max},
\end{align*}
\red{where $\{\omega_{k}^\bsB\}^{M}_{k=1}$ is the parametrize vector of the sparse representation \eqref{eq:b} for $\bsB$.}
\end{corollary}
}
\begin{remark}
When $\bsB$ is not sparse and just a general unstructured matrix, from \eqref{eq:rm} and Corollary \ref{th:important-express-cond_f}, it yields
\begin{equation}\label{eq:un 46}
	\begin{split}
		\red{\cond_{|\Omega_{QS}|,|\bsB|}(\bsA(\Omega_{QS}),\bsB )} &= \frac{1}{\|\bsX\|_{\max}}\bigg\| \left | \bsA^{-1}  \right| |\bsB| +|\bsA^{-1}||\bsA_{\bsD}||\bsX|
|\bsA^{-1}\bsA_{\bsL}||\bsX|\\
&\quad +|\bsA^{-1}||\bsA_{\bsL}\bsX|+|\bsA^{-1}||\bsA_{\bsU}\bsX|+|\bsA^{-1}\bsA_{\bsU}||\bsX|\\
&\quad +
\sum^{n-1}_{i=2}\Big|\bsA^{-1}
\begin{bmatrix}0&0\\\bsA(i+1:n,1:i-1)&0\\
\end{bmatrix}\bsX\Big|\\
&\quad+\sum^{n-1}_{i=2}\left|\bsA^{-1}
\begin{bmatrix}0&\bsA(1:i-1,i+1:n)\\0&0\\
\end{bmatrix}\bsX\right|  \bigg\|_{\max}.
	\end{split}
\end{equation}

\end{remark}

From Corollary \ref{th:important-express-cond_f}, we can claim  that $\cond_{|\Omega_{QS}|,|\Omega_{\bsB}|}(\bsA(\Omega_{QS}),\bsB(\Omega_\bsB) ) $ only depends on the entries of $\bsA$, $\bsX$, $\bsS_k$ and \red{$|\Omega_\bsB|$}.


\begin{proposition}\label{pro:scal-cond_f}
Let $\bsA\bsX=\bsB$, where $0\neq \bsX$ $\in \bsR^{n\times m}$, $\bsA$ $\in \R^{n\times m}$ is a \{1;1\}-quasiseparable \red{nonsingular} matrix, and $\bsB\in \R^{n \times m}$ is given by \eqref{eq:b}. Then, for any two vectors $\Omega_{QS}$ and $\Omega'_{QS}$ of quasiseparable parameters of $\bsA$, we have
     $$
    \red{\cond_{|\Omega_{QS}|,|\Omega_\bsB|}(\bsA(\Omega_{QS}),\bsB(\Omega_\bsB) )=\cond_{|\Omega^{'}_{QS}|,|\Omega_\bsB|}(\bsA(\Omega'_{QS}),\bsB(\Omega_\bsB) )}.
     $$
\end{proposition}

\end{remark}

\subsection{The Givens vector representation via tangents }\label{section:5.2}



The Givens-vector representation, introduced first in \cite{vandebril2005a}, is another important representation for quasiseparable matrices.\,\,Later, its minor variant called tangent-Givens-vector representation \eqref{eq:GV rep} was presented in \cite{Dopico2015Structured}. In this subsection, an explicit expression of the componentwise condition number for the solution of \eqref{eq:mul} with a quasiseparable matrix of coefficients with respect to the Givens-vector representation via tangents \eqref{eq:GV rep}  is presented in Theorem \ref{th:explict-expression-Gv}. The proof is straightforward from Theorem \ref{th:important- expres-compon}, Lemma \ref{lemm:derivati-partial-Givens}, Remark \ref{re:dve}. Therefore, we omit the proof of Theorem \ref{th:explict-expression-Gv}.

\red{Suppose that $\bsA$ has a  Givens-vector representation via tangents  \eqref{eq:GV rep} and $\bsB \in \R^{n\times m}$ \red{has the sparse representation} \eqref{eq:b}. Considering Definition \ref{def:com-paramter-condi-num}, we introduce  the weight vector
 \begin{equation}\label{eq:gv}
 \bsE_{GV}=\Big(\{\bse_{\bsl_{i}}\}^{n-1}_{i=2}, \{\bse_{\bsv_{i}}\}^{n-1}_{i=1}, \{\bse_{\bsd_{i}}\}^{n}_{i=1}, \{\bse_{\red{\bsw_{i}}}\}^{n-1}_{i=1}, \{\bse_{\bsu_{i}}\}^{n-1}_{i=2}\Big).
 \end{equation}
 with respect to the the Givens-vector representation via tangents  \eqref{eq:GV rep}  of $\bsA$. Therefore, using Definition \ref{def:com-paramter-condi-num}, we introduce the structure componentwise condition number for \eqref{eq:mul} with respect to the  Givens-vector representation via tangents  \eqref{eq:GV rep}  of $\bsA$ and the sparse representation \eqref{eq:b} of $\bsB$ in Theorem \ref{th:explict-expression-Gv}.
}

\begin{theorem}\label{th:explict-expression-Gv}
Let $\bsA \bsX = \bsB$, where $ \bsX$ $\in \R^{n \times m}$, \red{$\bsB$ $\in \R^{n \times m}$ has the sparse representation  \eqref{eq:b}} and $\bsA$ $\in \R^{n \times n}$ is a \{1;1\}-quasiseparable nonsingular matrix with a tangent-Givens-vector representation $\Omega_{GV}$ \eqref{eq:GV rep}, $\bsA=\bsA_{\bsL}+\bsA_{\bsD}+\bsA_{\bsU}$, with $\bsA_{\bsL}$ strictly lower triangular, $\bsA_{\bsD}$ diagonal, and $\bsA_{\bsU}$ strictly upper triangular. Let $\mathbf{0}\leq \bsf  \in \R^{M}$ and $\mathbf{0}\leq \bsE_{GV}$ $\in \R^{5n-6}$, \red{where  $\bsE_{GV}$ is given by \eqref{eq:gv}.} Suppose $\Omega_{GV} \cap \Omega_\bsB= \emptyset$, then
\begin{align*}
&\cond_{\bsE_{GV},\bsf}(\bsA(\Omega_{GV},\bsB(\Omega_{\bsB } ) )=\frac{1}{\|\bsX\|_{\max}}\Bigg\|  \sum^{M}_{k=1} \left | \bsA^{-1}  \bsS_k \right| \bsf_{k}  +  |\bsA^{-1}||\bsQ_{\bsd}||\bsX| +
|\bsA^{-1}\bsA_{\bsL}||\bsQ_{\bsv}||\bsX| \\
&\quad +  |\bsA^{-1}||\bsQ_{\bsw}||\bsA_{\bsU}\bsX| +
\sum^{n-1}_{i=2}\left| \bsA^{-1}\begin{bmatrix}0&0\\-\bss^{2}_{i}\bsA(i,1:i-1)&0\\\bsc^{2}_{i}\bsA(i+1:n,1:i-1)&0\\\end{bmatrix}\bsX\right| \left|\frac{\bse_{\bsl_{i}}}{\bsl_{i}}\right|\\
&\quad +\sum^{n-1}_{i=2}\left| \bsA^{-1} \begin{bmatrix} 0&-\bst^{2}_{i}\bsA(1:i-1,i)&\bsr^{2}_{i}\bsA(1:i-1,i+1:n)\\0&0&0\\ \end{bmatrix} \bsX\right| \left|\frac{\bse_{\bsu_{i}}}{\bsu_{i}}\right|  \Bigg\|_{\max},
\end{align*}
where
\begin{align}\label{eq:th5.1}
    \Omega_{GV}&=\Big(\{\bsl_{i}\}^{n-1}_{i=2}, \{\bsv_{i}\}^{n-1}_{i=1}, \{\bsd_{i}\}^{n}_{i=1}, \{\bsw_{i}\}^{n-1}_{i=1}, \{\bsu_{i}\}^{n-1}_{i=2}\Big),\\
   \bsQ_{\bsd}&=\diag\Big(\bse_{\bsd_{1}}, \ldots, \bse_{\bsd_{n}}\Big),   \bsQ_{\bsv}=\diag\Big(\frac{\bse_{\bsv_{1}}}{\bsv_{1}}, \ldots, \frac{\bse_{\bsv_{n-1}}}{\bsv_{n-1}},1\Big), \bsQ_{\red{\bsw}}=\diag\Big(\frac{\bse_{\red{\bsw_{1}}}}{\red{\bsw_{1}}}, \ldots, \frac{\bse_{\red{\bsw_{n-1}}}}{\red{\bsw_{n-1}}}, 1\Big).\notag
\end{align}
\end{theorem}

\begin{remark}\label{remark:2}
As in Remark \ref{remark:1}, we consider the situation that multiple right-hand sides $\bsB$ is not sparse \red{but} general unstructured. According to Theorem \ref{th:explict-expression-Gv}, the explicit expression for $\cond_{\bsE_{GV},\bsF}(\bsA(\Omega_{GV}),\bsB)$ of \red{the multiple right-hand sides linear system} \eqref{eq:mul}, where $\bsA$ is a nonsingular \{1;1\}-quasiseparable matrix with a quasiseparable representation $\Omega_{GV}$ and $\bsB$ is a general dense matrix, can be characterized by
\begin{align}\notag 
\cond_{\bsE_{GV},\bsF}&(\bsA(\Omega_{GV}),\bsB ):=\lim_{\eta \rightarrow 0}\sup\Bigg  \{ \frac{1}{\eta}\frac{\|\delta \bsX\|_{\max  }}{\|\bsX\|_{\max  }}:\bsA(\Omega_{GV} +\delta\Omega_{GV} )(\bsX+\delta\bsX)=\bsB+\delta\bsB,  \nonumber \\
&\hspace{4.2cm}  |\delta \Omega_{GV} |\leq \eta \bsE_{GV}, |\delta \bsB |\leq \eta\bsF
\Bigg \} \nonumber \\
&=\frac{1}{\|\bsX\|_{\max}}\Bigg\|  |\bsA^{-1}|\bsF   +  |\bsA^{-1}||\bsQ_{\bsd}||\bsX| +
|\bsA^{-1}\bsA_{\bsL}||\bsQ_{\bsv}||\bsX|  \notag \\
&\, +  |\bsA^{-1}||\bsQ_{\bsw}||\bsA_{\bsU}\bsX| +
\sum^{n-1}_{i=2}\left| \bsA^{-1}\begin{bmatrix}0&0\\-\bss^{2}_{i}\bsA(i,1:i-1)&0\\\bsc^{2}_{i}\bsA(i+1:n,1:i-1)&0\\\end{bmatrix}\bsX\right| \left|\frac{\bse_{\bsl_{i}}}{\bsl_{i}}\right|\nonumber \\
&\, +\sum^{n-1}_{i=2}\left| \bsA^{-1} \begin{bmatrix} 0&-\bst^{2}_{i}\bsA(1:i-1,i)&\bsr^{2}_{i}\bsA(1:i-1,i+1:n)\\0&0&0\\ \end{bmatrix} \bsX\right| \left|\frac{\bse_{\bsu_{i}}}{\bsu_{i}}\right|  \Bigg\|_{\max}. \nonumber
\end{align}
For the linear system $\bsA \bsx =\bsb $ with a \{1;1\}-quasiseparable coefficient matrix $\bsA(\Omega_{GV})$ in the Givens-vector representation via tangents $\Omega_{GV}$, Dopico and Pom{\'e}s \cite{Dopico2016} introduced the componentwise condition number of $\bsx$ with respect to $\Omega_{GV}$  and the corresponding explicit formula \cite[Theorem 5.6]{Dopico2016} is given below
\begin{align*}
&{\rm cond}_{\bsE_{GV},\bsy  }(\bsA(\Omega_{GV},\bsb )=\frac{1}{\|\bsx\|_{\infty}}\Bigg\|  |\bsA^{-1}|\bsy    +  |\bsA^{-1}||\bsQ_{\bsd}||\bsx|+
|\bsA^{-1}\bsA_{\bsL}||\bsQ_{\bsv}||\bsx| \\
&\quad + |\bsA^{-1}||\bsQ_{\bsw}||\bsA_{\bsU}\bsx| +
\sum^{n-1}_{i=2}\left| \bsA^{-1}\begin{bmatrix}0&0\\-\bss^{2}_{i}\bsA(i,1:i-1)&0\\\bsc^{2}_{i}\bsA(i+1:n,1:i-1)&0\\\end{bmatrix}\bsx\right| \left|\frac{\bse_{\bsl_{i}}}{\bsl_{i}}\right|\\
&\quad + \sum^{n-1}_{i=\red{2}}\left| \bsA^{-1} \begin{bmatrix} 0&-\bst^{2}_{\red{i}}\bsA(1:i-1,i)&\bsr^{2}_{i}\bsA(1:i-1,i+1:n)\\0&0&0\\ \end{bmatrix} \bsx\right| \left|\frac{\bse_{\bsu_{i}}}{\bsu_{i}}\right|  \Bigg\|_{\infty },
\end{align*}
where $  \Omega_{GV}$, $\bsE_{GV}$, $ \bsQ_{\bsd}$, $\bsQ_{\bsv}$ and $\bsQ_{\bsw}$ is defined in \eqref{eq:th5.1}. For the case $m=1$, when $\bsB=\bsb$ and $\bsF=\bsy $, it is not difficult to see that the condition number of the multi-right side linear system $\cond_{\bsE_{GV},\bsF}(\bsA(\Omega_{GV}),\bsB )$ is identical to ${\rm cond}_{\bsE_{GV},\bsy }(\bsA(\Omega_{GV}),\bsb )$.
\end{remark}

\red{In view of Theorem \ref{th:important- expres-compon}, a natural  choice of $\bsE_{GV}$ \eqref{eq:gv} is $\bsE_{GV} =|\Omega_{GV}|$, where  $\Omega_{GV}$ is given by \eqref{eq:GV rep} and is the Givens-vector representation via tangents introduced in Definition \ref{def:givens } of $\bsA$}. Under this situation, noting $\bsQ_{\bsd}=|\bsA_{\bsD}|$, $\bsQ_{\bsv}=\bsQ_{\bse}=\bsI$, and $\big|\frac{\bse_{\bsl_{i}}}{\bsl_{i}}\big|=\big|\frac{\bse_{\bsu_{i}}}{\bsu_{i}}\big|=1$, we have
\red{
\begin{corollary}\label{th:explic-cond_f(GV)}
Let $\bsA \bsX=\bsB$, where $\mathbf{0}\neq \bsX$ $\in\R^{n \times m}$, $\bsB$ $\in \R^{n \times m}$, and $\bsA$ $\in\R^{n \times n}$ is a \{1;1\}-quasiseparable \red{nonsingular} matrix with a tangent-Givens-vector representation $\Omega_{GV}$ \eqref{eq:GV rep}, $\bsB \in \R^{n\times m}$ \red{has the sparse representation} \eqref{eq:b}, and $\bsA=\bsA_{\bsL}+\bsA_{\bsD}+\bsA_{\bsU}$, with $\bsA_{\bsL}$ strictly lower triangular, $\bsA_{\bsD}$ diagonal and $\bsA_{\bsU}$ strictly upper triangular. Suppose $\Omega_{GV} \cap\, \Omega_\bsB= \emptyset$\red{.} \red{Then}
\begin{align*}
\red{\cond_{|\Omega_{GV}|,|\Omega_\bsB|} }&(\bsA(\Omega_{GV}),\bsB(\Omega_{\bsB }))=\frac{1}{\|\bsX\|_{\max}}\Bigg\|  \sum^{M}_{k=1} \left | \bsA^{-1}  \bsS_k \right| \red{|\bsw^{\bsB}_{k}|}   +  |\bsA^{-1}||\bsA_{\bsD}||\bsX| +
|\bsA^{-1}\bsA_{\bsL}||\bsX| \\
&\, + |\bsA^{-1}||\bsA_{\bsU}\bsX|  +
\sum^{n-1}_{i=2}\left| \bsA^{-1}\begin{bmatrix}0&0\\-\bss^{2}_{i}\bsA(i,1:i-1)&0\\\bsc^{2}_{i}\bsA(i+1:n,1:i-1)&0\\\end{bmatrix}\bsX\right|\\
&\, +\sum^{n-1}_{i=1}\left| \bsA^{-1} \begin{bmatrix} 0&-\bst^{2}_{\red{i}}\bsA(1:i-1,i)&\bsr^{2}_{i}\bsA(1:i-1,i+1:n)\\0&0&0\\ \end{bmatrix} \bsX\right| \Bigg\|_{\max}.
\end{align*}
\end{corollary}
}

\begin{remark}
	When $\bsB$ is not sparse and just a general unstructured matrix, it follows from \eqref{eq:rm} and Corollary \ref{th:important-express-cond_f} that
	\begin{equation}\label{eq:410}
		\begin{split}
			\red{\cond_{|\Omega_{GV}|,|\bsB|} }&(\bsA(\Omega_{GV}),\bsB )=\frac{1}{\|\bsX\|_{\max}}\Bigg\|   \left | \bsA^{-1}   \right| |\bsB|   +  |\bsA^{-1}||\bsA_{\bsD}||\bsX| +
|\bsA^{-1}\bsA_{\bsL}||\bsX| \\
&\, + |\bsA^{-1}||\bsA_{\bsU}\bsX|  +
\sum^{n-1}_{i=2}\left| \bsA^{-1}\begin{bmatrix}0&0\\-\bss^{2}_{i}\bsA(i,1:i-1)&0\\\bsc^{2}_{i}\bsA(i+1:n,1:i-1)&0\\\end{bmatrix}\bsX\right|\\
&\, +\sum^{n-1}_{i=1}\left| \bsA^{-1} \begin{bmatrix} 0&-\bst^{2}_{\red{i}}\bsA(1:i-1,i)&\bsr^{2}_{i}\bsA(1:i-1,i+1:n)\\0&0&0\\ \end{bmatrix} \bsX\right| \Bigg\|_{\max}.
		\end{split}
	\end{equation}
\end{remark}

Comparing the expressions in Corollaries  \ref{th:explic-cond_f(GV)} and \ref{th:important-express-cond_f}, we know that the structured condition number \red{$\cond_{|\Omega_{GV}|,|\Omega_\bsB|}(\bsA(\Omega_{GV}),\bsB(\Omega_{\bsB }))$} \red{with respect to} the Givens-vector representation \eqref{eq:GV rep} depends on not only all elements of $\bsA$ but also  the parameters $\{\bsc_{i},\bss_{i}\}$ and $\{\bsr_{i},\bst_{i}\}$. On the contrary, the structured condition number \red{$\cond_{|\Omega_{QS}|,|\Omega_\bsB|}\big(\bsA(\Omega_{QS}),\bsB(\Omega_{\bsB })\big)$}  \red{with respect to the quasiseparable representation \eqref{eq:qs general} does not rely on these quasiseparable representation parameters defining $A$}.

\section{Relationships between different condition numbers for multiple right-hand side  linear systems with parametrized coefficient matrices }\label{section:6}

In this section, in Proposition \ref{pro:the comparison-cond_f and cond_A,f}
we will compare $\cond_{|\Omega_{QS}|,|\Omega_\bsB|}(\bsA(\Omega_{QS}),\bsB(\Omega_\bsB))$ given in Corollary \ref{th:important-express-cond_f}   with the unstructured componentwise condition number \red{$\cond_{|\bsA|,|\Omega_\bsB| }(\bsA,\bsB(\Omega_{\bsB }))$} given by \eqref{eq:cond_unstruWithSparseRH}. \red{The relationship between 
$$\cond_{|\Omega_{QS}|,|\Omega_\bsB|}(\bsA(\Omega_{QS}),\bsB(\Omega_\bsB) ) \mbox{ and } \cond_{|\Omega_{GV}|,|\Omega_\bsB| }(\bsA(\Omega_{GV}),\bsB(\Omega_{\bsB }))
$$ given in Corollary \ref{th:explic-cond_f(GV)} is also investigated in Theorem \ref{th:the inequ-cond{QS} and condeff{QS}}. The effective condition number $\cond_{\rm eff}(\bsA(\Omega_{QS}),\bsB(\Omega_\bsB) )$ is introduced in Definition \ref{def:cond_eff}, which can bound the exact condition numbers $\cond_{|\Omega_{QS}|,|\Omega_\bsB|}(\bsA(\Omega_{QS}),\bsB(\Omega_\bsB) )$ and $\cond_{|\Omega_{GV}|,|\Omega_\bsB| }(\bsA(\Omega_{GV}),\bsB(\Omega_{\bsB }))$ up to the order of $n$.}  An important feature of \red{the effective condition number} is that it can be computed in $\Oh(mn)$ operations based on some previous results on \{1;1\}-quasiseparable matrices from \cite{eidelman1999on} and \cite{eidelman1999linear}. 

In Proposition \ref{pro:the comparison-cond_f and cond_A,f}, we will investigate the relationship between the structured componentwise condition number \red{$\cond_{|\Omega_{QS}|,|\Omega_\bsB|}(\bsA(\Omega_{QS}),\bsB(\Omega_\bsB))$} and the corresponding unstructured one \red{$\cond_{|\bsA|,|\Omega_\bsB|}(\bsA,\bsB(\Omega_\bsB) )$} given in Corollary \ref{co:3.1}. 


%

\begin{proposition}\label{pro:the comparison-cond_f and cond_A,f}
Let $\bsA \bsX=\bsB$, where $\bsA$ $\in \R^{n\times n}$ is a \{1;1\}-quasiseparable nonsingular matrix, $\bsB$ $\in \R^{n\times m}$ has the sparse representation as \eqref{eq:b}, and $0\neq \bsX$ $\in \R^{n\times m}$. Let $\Omega_{QS}$ be a quasiseparable representation of $\bsA$ and $\Omega_\bsB $ be a sparse representation of $\bsB $. Then the following relation holds,
$$
\red{\cond_{|\Omega_{QS}|,|\Omega_\bsB|}(\bsA(\Omega_{QS}),\bsB(\Omega_{\bsB }))}\leq n\,\, \cond_{|\bsA|,\red{|\Omega_\bsB|} }(\bsA,\bsB(\Omega_{\bsB }) ).
$$
\end{proposition}

\proof From Corollary  \ref{th:important-express-cond_f} and using standard properties of absolute values and norms we obtain:
\begin{align*}
\cond_{|\Omega_{QS}|,\red{|\Omega_\bsB|}}(\bsA(\Omega_{QS}),\bsB(\Omega_{\bsB })) &\leq \frac{1}{\|\bsX\|_{\max}}\bigg\|\sum^{M}_{k=1} \left | \bsA^{-1}  \bsS_k \right| \red{|\bsw^{\bsB}_{k}|}  + |\bsA^{-1}||\bsA_{\bsD}||\bsX|+ 2|\bsA^{-1}||\bsA_{\bsL}||\bsX|  \\
&\quad +2|\bsA^{-1}||\bsA_{\bsU}||\bsX| +\sum^{n-1}_{i=2}|\bsA^{-1}||\bsA_{\bsL}||\bsX|+\sum^{n-1}_{i=2}|\bsA^{-1}||\bsA_{\bsU}||\bsX| \bigg\|_{\max}\\
&\quad \leq\frac{n}{\|\bsX\|_{\max}}\Big\| \sum^{M}_{k=1} \left | \bsA^{-1}  \bsS_k \right| \red{|\bsw^{\bsB}_{k}|}  +  |\bsA^{-1}||\bsA||\bsX| \Big\|_{\max},
\end{align*}
which finishes the proof.
\eproof

From the above  proposition, the structured condition number $\cond_{|\Omega_{QS}|,\red{|\Omega_\bsB|}}(\bsA(\Omega_{QS}),\bsB(\Omega_{\bsB }))  $ is smaller than the unstructured condition number $ \cond_{|\bsA|,\red{|\Omega_\bsB|} }(\bsA,\bsB(\Omega_{\bsB }) )$, except for a factor $n$. In addition, we will see in the numerical experiments presented in Section \ref{section:8} that it can be much smaller.

 In the following theorem, we will study the relationship between $\cond_{|\Omega_{GV}|,\red{|\Omega_\bsB|}}(\bsA(\Omega_{GV}),\bsB(\Omega_\bsB))$ and $\cond_{|\Omega_{QS}|,\red{|\Omega_\bsB|}}(\bsA(\Omega_{QS}),\bsB(\Omega_\bsB) )$. The similar conclusions had been obtained for eigenvalue, generalized eigenvalue computations and linear system solving  for \{1;1\}-quasiseparable matrices in \cite{Dopico2015Structured,Dopico2016,diaomeng}, respectively. Before that, we need \red{to review} Lemma \ref{lem:6ineq}, which describes the perturbation magnitude relationship between the Givens-vector representation via tangents given in Definition \ref{def:givens } and the quasiseparable representation given in Definition \ref{th:{1;1}}. \red{In the rest of this paper,   we will adopt the following brevity notations}:
\red{
\begin{align}\label{notaion:1}
\cond^{QS}_{|\bsB|} := \cond_{|\Omega_{QS}|,\,|\bsB|}(\bsA(\Omega_{QS}),|\bsB|),\quad \cond^{QS}_{|\Omega_{\bsB}|} := \cond_{|\Omega_{QS}|,|\Omega_{\bsB}|}(\bsA(\Omega_{QS}),|\Omega_{\bsB}|),
\end{align}
where $\cond_{|\Omega_{QS}|,\,|\bsB|}(\bsA(\Omega_{QS}),|\bsB|)$ and $\cond_{|\Omega_{QS}|,|\Omega_{\bsB}|}(\bsA(\Omega_{QS}),|\Omega_{\bsB}|)$ are given in  \eqref{eq:un 46} and Corollary \ref{th:important-express-cond_f}, respectively;}
\red{
\begin{align}\label{notaion:2}
\cond^{GV}_{|\bsB|} := \cond_{|\Omega_{GV}|,\,|\bsB|}(\bsA(\Omega_{GV}),|\bsB|),\quad \cond^{GV}_{|\Omega_{\bsB}|} := \cond_{|\Omega_{GV}|,|\Omega_{\bsB}|}(\bsA(\Omega_{GV}),|\Omega_{\bsB}|),
\end{align}
where $\cond_{|\Omega_{GV}|,\,|\bsB|}(\bsA(\Omega_{GV}),|\bsB|)$ and $\cond_{|\Omega_{GV}|,|\Omega_{\bsB}|}(\bsA(\Omega_{GV}),|\Omega_{\bsB}|)$ are given in  \eqref{eq:410} and Corollary \ref{th:explic-cond_f(GV)}, respectively;}
\red{
\begin{align}\label{notaion:3}
\cond_{|\bsB|}:= \cond_{|\bsA|,|\bsB|}(\bsA,\bsB),\quad \quad \cond_{|\Omega_{\bsB}|}:=\cond_{|\bsA|,|\Omega_{\bsB}|}(\bsA,\Omega_{\bsB}),
\end{align}
where $ \cond_{|\bsA|,|\bsB|}(\bsA,\bsB)$ and $\cond_{|\bsA|,|\Omega_{\bsB}|}(\bsA,\Omega_{\bsB})$ are given in  \eqref{th:1the expression-compon} and \eqref{eq:cond_unstruWithSparseRH}, respectively. }

 \begin{lemma}\cite[Lemma 6.2]{Dopico2015Structured}\label{lem:6ineq}
	With the notations before, we have
	$$
	\left|\delta \Omega^{GV}_{QS} \right |\leq \eta \left  |\Omega^{GV}_{QS} \right | \Rightarrow |\delta{'}\Omega_{GV})|\leq (3(n-2)\eta+\Oh(\eta^{2})) \left |\Omega_{GV}\right |.
	$$
\end{lemma}

\begin{theorem}\label{th:cond(GV)and cond(QS)}
Let $\bsA\bsX=\bsB$, where $ \bsX$ $\in \R^{n \times m}$\red{,} $\bsA$ $\in\R^{n\times n}$ is a nonsingular \{1;1\}-quasiseparable matrix \red{with a tangent-Givens-vector representation $\Omega_{GV}$ \eqref{eq:GV rep} and a quasiseparable representation $\Omega_{QS}$} \eqref{eq:qs general}. Suppose $\bsB$ \red{has the sparse representation} \eqref{eq:b}, \red{ and assume that $\Omega_{QS} \cap \Omega_\bsB= \emptyset$ and $\Omega_{GV} \cap \Omega_\bsB= \emptyset$}. Then the following relationship holds:
$$ \red{\cond^{GV}_{|\Omega_\bsB|}\,\,\leq \,\,\cond^{QS}_{|\Omega_\bsB|} \leq \,\, (3n-2)\, \cond^{GV}_{|\Omega_\bsB|}}.
$$
\end{theorem}

\proof
First, we show \red{$\cond^{GV}_{|\Omega_\bsB|}\,\,\leq \,\,\cond^{QS}_{|\Omega_\bsB|}$}.  In view of the expressions of \red{$\cond^{GV}_{|\Omega_\bsB|}$} and \red{$\cond^{QS}_{|\Omega_\bsB|}$} given by Corollaries \ref{th:explic-cond_f(GV)} and \ref{th:important-express-cond_f}, respectively, we only \red{need to} compare the last two terms in the expressions of \red{$\cond^{GV}_{|\Omega_\bsB|}$} with the correspond parts of \red{$\cond^{QS}_{|\Omega_\bsB|}$} as follows:
\begin{align*}
	\sum^{n-1}_{i=2}&\left| \bsA^{-1}\begin{bmatrix}0&0\\-\bss^{2}_{i}\bsA(i,1:i-1)&0\\\bsc^{2}_{i}\bsA(i+1:n,1:i-1)&0\\\end{bmatrix}\bsX\right|
	\leq
	\sum^{n-1}_{i=2}\left| \bsA^{-1}\begin{bmatrix}0&0\\  \bsA(i,1:i-1)&0\\0&0\\\end{bmatrix}\bsX\right|\\
	&\quad +\sum^{n-1}_{i=2}\left| \bsA^{-1}\begin{bmatrix}0&0\\0 &0\\\ \bsA(i+1:n,1:i-1)&0\\\end{bmatrix}\bsX\right|\\
	&\quad = \sum^{n-1}_{i=2}\left| \bsA^{-1} e_i \cdot  \bsA_\bsL(i,: ) \bsX\right| +\sum^{n-1}_{i=2}\left| \bsA^{-1}\begin{bmatrix}0&0 \\ \bsA(i+1:n,1:i-1)&0\\\end{bmatrix}\bsX\right|\\
	&\quad \leq  \left| \bsA^{-1} \right | ~ \left | \bsA_\bsL \bsX\right| +\sum^{n-1}_{i=2}\left| \bsA^{-1}\begin{bmatrix}0&0 \\ \bsA(i+1:n,1:i-1)&0\\\end{bmatrix}\bsX\right|,
\end{align*}
where in the last two inequality we use \eqref{eq:out}. \red{Similarly}, we can prove that
\begin{align*}
	\sum^{n-1}_{i=1}&\left| \bsA^{-1} \begin{bmatrix} 0&-\bst^{2}_{i}\bsA(1:i-1,i)&\bsr^{2}_{i}\bsA(1:i-1,i+1:n)\\0&0&0\\ \end{bmatrix} \bsX\right|
	\\
&\leq \quad  |\bsA^{-1}\bsA_{\bsU}||\bsX| +\sum^{n-1}_{i=2}\left|\bsA^{-1}
\begin{bmatrix}0&\bsA(1:i-1,i+1:n)\\0&0\\
\end{bmatrix}\bsX\right|.
\end{align*}
Thus we complete the first part proof.  Now we turn to the second part of the proof.  Note that from Definition \ref{def:com-paramter-condi-num} and from Lemma \ref{lem:6ineq} we have
\begin{align*}
\red{\cond^{QS}_{|\Omega_\bsB|}}&\leq \lim_{\eta \rightarrow 0}
\sup \Bigg\{ \frac{\|\delta\bsX\|_{\max}}{\eta\|\bsX\|_{\max}}:(\bsA(\Omega_{GV}+\delta\Omega_{GV}))(\bsX+\delta\bsX)=\bsB(\Omega_{\bsB}+\delta \Omega_{\bsB} ), \\
&|\delta\Omega_{GV}|\leq \left  [ 3(n-2)\,\eta+\Oh(\eta^{2}) \right ] \,|\Omega_{GV}|, |\delta\,\bsB|\leq \left  [ 3(n-2)\,\eta+\Oh(\eta^{2} ) \right ] \,\red{|\Omega_\bsB|}
\Bigg\}.
\end{align*}
By considering the change of variable $\eta' = (3(n-2)\eta + \Oh(\eta^{2}))$, we obtain
\begin{align*}
\red{\cond^{QS}_{|\Omega_\bsB|}\leq}\lim_{\eta'\rightarrow 0}&\Bigg\{ \frac{\|\delta\bsX\|_{\max}}{\eta\|\bsX\|_{\max}}:(\bsA(\Omega_{GV}+\delta\,\Omega_{GV}))(\bsX+\delta\,\bsX)=\bsB(\Omega_{\bsB}+\delta \Omega_{\bsB} ) ,\\
&\quad \quad \quad \quad \quad\quad |\delta\,\Omega_{GV}|\leq \eta'\,|\Omega_{GV}|, |\delta\,\Omega_{\bsB}|\leq \eta'\,\red{|\Omega_{\bsB}|}
\Bigg\}\\
&=3(n-2)\,\,\red{\cond^{GV}_{|\Omega_\bsB|}}.
\hspace{4.6cm} \Box
\end{align*}

Although the explicit expressions of \red{$\cond^{QS}_{|\Omega_\bsB|}$} and \red{$\cond^{GV}_{|\Omega_\bsB|}$} have be derived in Corollaries \ref{th:important-express-cond_f} and \ref{th:explic-cond_f(GV)}, respectively, the formulas involve two sums, which can result in expensive computational costs. Therefore it is desirable  to consider other condition numbers having the similar contributions of \red{$\cond^{QS}_{|\Omega_\bsB|}$} or \red{$\cond^{GV}_{|\Omega_\bsB|}$}. In Definition \ref{def:cond_eff} we will propose the effective condition number $\cond_{\rm eff}(\bsA(\Omega_{QS}), \bsB(\Omega_{\bsB}))$. Moreover, we will investigate the relationship between $\cond_{\rm eff}(\bsA(\Omega_{QS}), \bsB(\Omega_{\bsB }) )$ and \red{$\cond^{QS}_{|\Omega_\bsB|}$/$\cond^{GV}_{|\Omega_\bsB|}$} in Theorem \ref{th:the inequ-cond{QS} and condeff{QS}} and Theorem \ref{th:cond{GV}-cond{QS}-cond_eff}.

\begin{definition}\label{def:cond_eff}
Let $\bsA\bsX=\bsB$, where $\mathbf{0} \neq \bsX$ $\in \R^{n \times m}$, \red{$\bsB$ has the sparse representation  \eqref{eq:b}}, and $\bsA$ $\in \R^{n\times n}$ is a nonsingular \{1;1\}-quasiseparable matrix \red{with a quasiseparable representation $\Omega_{QS}$} \eqref{eq:qs general} such that $\bsA=\bsA_{\bsL}+\bsA_{\bsD}+\bsA_{\bsU}$, with $\bsA_{\bsL}$ strictly lower triangular, $\bsA_{\bsD}$ diagonal, and $\bsA_{\bsU}$ strictly upper triangular.\,\red{Assume that $\Omega_{QS} \cap \,\Omega_\bsB= \emptyset$}. Then, we define the effective relative condition number $\cond_{\rm eff}(\bsA(\Omega_{QS}), \bsB(\Omega_{\bsB }) )$ for the solution of $\bsA \bsX=\bsB$ as
\begin{align}\notag 
\cond_{\rm eff}(\bsA(\Omega_{QS}), \bsB(\Omega_{\bsB }) ):=&\frac{1}{\|\bsX\|_{\max}}   \Bigg\|
 \sum^{M}_{k=1} \left | \bsA^{-1}  \bsS_k \right| \red{|\bsw^{\bsB}_{k}|}  + |\bsA^{-1}||\bsA_{\bsD}||\bsX| + |\bsA^{-1}||\bsA_{\bsL}\bsX| \nonumber\\ &+ |\bsA^{-1}\bsA_{\bsL}||\bsX| + |\bsA^{-1}||\bsA_{\bsU}\bsX| + |\bsA^{-1}\bsA_{\bsU}||\bsX|            \Bigg\|_{\max}. \notag
\end{align}
\end{definition}

\red{
\begin{remark}
In \cite[Algorithm 5.2]{eidelman1999linear}, it was shown that the inverse of a nonsingular quasiseparable matrix $\bsA$ $\in \R^{n\times n}$ is also quasiseparable and its quasiseparable representation can be obtained in $58(n-2)+20$ flops. Clearly, $\bsA_{\bsL}$ and $\bsA_{\bsU}$ are quasiseparable. Furthermore, from \cite[Algorithm 4.4]{eidelman1999on} the matrix-vector multiplication $\bsA \bsv$ can be calculated in $\Oh(n)$, where $\bsA \in \R^{n \times n}$ is a quasiseparable matrix. Therefore, $\cond_{\rm eff}(\bsA(\Omega_{QS}), \bsB(\Omega_{\bsB }) )$ can be computed in $\Oh(mn)$ flops. The detailed descriptions of the routine for computing $\cond_{\rm eff}(\bsA(\Omega_{QS}), \bsB(\Omega_{\bsB }) )$  are omitted.
\end{remark}
}

In the following theorem, we will investigate the relationship between \red{$\cond_{\rm eff}(\bsA(\Omega_{QS}), \bsB(\Omega_{\bsB }))$ and $\cond^{QS}_{|\Omega_\bsB|}$.}

\begin{theorem}\label{th:the inequ-cond{QS} and condeff{QS}}
Let $\bsA \bsX=\bsB$, where $\mathbf{0} \neq \bsX \in \R^{n \times m}$ and $\bsA \in \R^{n \times n}$ is a nonsingular \{1;1\}-quasiseparable matrix \red{with a quasiseparable representation $\Omega_{QS}$} \eqref{eq:qs general}, \red{$\bsB\in \R^{n\times m}$ has the sparse representation  \eqref{eq:b}. Let $\Omega_{QS} \cap \Omega_\bsB= \emptyset$}. Then the following relations hold:
$$
\cond_{\rm eff}(\bsA(\Omega_{QS}), \bsB(\Omega_{\bsB }) )
\leq \red{\cond^{QS}_{|\Omega_\bsB|}} \leq (n-1) \,\cond_{\rm eff}(\bsA(\Omega_{QS}), \bsB(\Omega_{\bsB }) ).
$$
\end{theorem}
\proof Clearly, \red{$\cond_{\rm eff}(\bsA(\Omega_{QS}), \bsB(\Omega_{\bsB }) )
\leq \cond^{QS}_{|\Omega_\bsB|}$ is easily obtained from Corollary \ref{th:important-express-cond_f}
  and Definition \ref{def:cond_eff}}. On the other hand, it yields
\begin{align*}
\bigg|\bsA^{-1}\begin{bmatrix}0&0\\\bsA(i+1:n,1:i-1)&0\\
              \end{bmatrix}\bsX \bigg|
              &=\bigg|\bsA^{-1}\begin{bmatrix}0&0\\\bsA_{\bsL}(i+1:n,1:i-1)&0\\
              \end{bmatrix}\bsX  \bigg|\\
              &=\bigg|\bsA^{-1}\begin{bmatrix}0&0\\ \bsA_{\bsL}(i+1:n,1:i-1)&\bsA_{\bsL}(i+1:n,i:n)
              \end{bmatrix} \bsX\\
              &\quad +\bsA^{-1}\begin{bmatrix}0&0\\0&-\bsA_{\bsL}(i+1:n,i:n)
              \end{bmatrix}\bsX\bigg|\\
              &\leq |\bsA^{-1}|\bigg|\begin{bmatrix}0 \\ \bsA_{\bsL}(i+1:n,:)
              \end{bmatrix}\bsX\bigg|\\
              &\quad +\bigg|\bsA^{-1}\begin{bmatrix}0&0\\0&\bsA_{\bsL}(i+1:n,i:n)\\
              \end{bmatrix}  \bigg||\bsX|\\
              &\leq \left |\bsA^{-1} \right | \left |\bsA_{\bsL} \bsX \right |+\left |\bsA^{-1}\bsA_{\bsL}\right | \left|\bsX \right|.
\end{align*}
Thus we can deduce that
\begin{equation}\label{eq:7-1}
\sum^{n-1}_{i=2}\bigg| \bsA^{-1}\begin{bmatrix} 0&0\\\bsA(i+1:n,1:i-1)&0
                               \end{bmatrix}\bsX \bigg| \leq (n-2)\, \left |\bsA^{-1} \right | \left|\bsA_{\bsL}\bsX \right | + (n-2)\, \left |\bsA^{-1}\bsA_{\bsL} \right| \left|\bsX \right|.
\end{equation}
Similarly, it can be derived that
\begin{equation}\label{eq:7-2}
\sum^{n-1}_{i=2}\left| \bsA^{-1}\begin{bmatrix} 0&\bsA^{-1}(1:i-1,i+1:n)\\0&0\\
                               \end{bmatrix}\bsX \right| \leq (n-2)\,\left |\bsA^{-1} \right | \left|\bsA_{\bsU}\bsX \right |+(n-2) \, \left|\bsA^{-1}\bsA_{\bsU} \right|
                               \left|\bsX \right |.
\end{equation}
Combing  (\ref{eq:7-1}) and (\ref{eq:7-2}), we finish the proof.
\eproof


From Theorem \ref{th:cond(GV)and cond(QS)} and Theorem \ref{th:the inequ-cond{QS} and condeff{QS}}, it is not difficult to deduce the relationship between $\cond_{\rm eff}(\bsA(\Omega_{QS}), \bsB(\Omega_{\bsB }))$ and \red{$\cond^{GV}_{|\Omega_\bsB|}$} in the following theorem.

\begin{theorem}\label{th:cond{GV}-cond{QS}-cond_eff}
Let $\bsA \bsX=\bsB$, where $\bsA$ $\in \R^{n \times n}$ is a nonsingular \{1;1\}-quasiseparable matrix with tangent-Givens-vector representation $\Omega_{GV}$ \eqref{eq:GV rep}, \red{$\bsB \in \R^{n\times m}$ has the sparse representation \eqref{eq:b}.\,\,For any quasiseparable representation $\Omega_{QS}$ of $\bsA$, \red{if} $\Omega_{QS} \cap \Omega_\bsB= \emptyset$ and $\Omega_{GV} \cap \,\Omega_\bsB= \emptyset$}, then the following relations hold:
\begin{align*}
\frac{ \cond_{\rm eff}(\bsA(\Omega_{QS}), \bsB(\Omega_{\bsB }) ) }{3(n-2)}\,\leq\,\red{\cond^{QS}_{|\Omega_\bsB|}} \,\leq\, (n-1) \, \cond_{\rm eff}(\bsA(\Omega_{QS}), \bsB(\Omega_{\bsB }) ) .
\end{align*}
\end{theorem}

In view of Theorem \ref{th:the inequ-cond{QS} and condeff{QS}} and Theorem \ref{th:cond{GV}-cond{QS}-cond_eff}, it is easy to see that the structured condition numbers \red{$\cond^{QS}_{|\Omega_\bsB|}$} and \red{$\cond^{GV}_{|\Omega_\bsB|}$} can be bounded by the effective condition number up to a factor of order $n$.

\section{Numerical experiments}\label{section:8}
\red{In this section, we do some numerical examples to illustrate the theoretical results for the multiple right-hand side linear system (\ref{eq:mul}). All the numerical experiments are carried out by {\sc Matlab} R2018a, with the machine epsilon $\mu\approx 2.2 \times 10^{-16}$. Recall that the condition number  notations defined in \eqref{notaion:1}, \eqref{notaion:2},  and \eqref{notaion:3}. For a given quasiseparable matrix $\bsA \in \R^{n \times n}$ and a multiple right-hand sides $\bsB \in \R^{m\times n}$, where $\bsB$ may be sparse or a general unstructured dense matrix, we compute the solution $\bsX$ to \eqref{eq:mul} by $\bsX=\bsA \backslash \bsB$ in  {\sc Matlab}. All condition numbers presented in this paper are computed directly from their explicit expressions in  {\sc Matlab}.
}

\begin{example}\label{example:1}
\red{Let
$$\bsA=\begin{bmatrix}1&-2.9442&0&0&0\\
                7.2688\cdot 10^4&1&1.4383\cdot 10 &0&0\\
                -2.6958\cdot 10^6&-3.0344\cdot 10^2&1&3.2519\cdot 10^{-1}&0\\
                -2.9947\cdot 10^9&-3.3709\cdot 10^5&2.9387\cdot 10^2&1&-7.5493\cdot 10^{-1}\\
                -8.5754\cdot 10^{12}&-9.6526\cdot 10^8&8.4150\cdot 10^5&-7.8728\cdot 10^2&1
\end{bmatrix},$$
$$\bsB_1 = \begin{bmatrix}
         1.0933 &0\\
         1.1093  &0\\
        -8.6365\cdot 10^{-1}&0\\
        0&7.7359\cdot 10^{-2}\\
          0&-1.2141
         \end{bmatrix},\,\,\,\,\,\,\,\,\,\,\,\,\,\,\,\,\,
\bsB_2=\begin{bmatrix}
            1.0000\cdot 10^{-3}&1\\
                   1&1\\
                   1&	1\\
                   1&	1\\
                   1&	1\\
         \end{bmatrix},
$$
where $\bsA$ is a \{1;1\}-quasiseparable matrix, $\bsB_1$ is a sparse matrix with a sparsity of 0.5 and $\bsB_2$ is a dense matrix.
In Table \ref{ta:nx0}, we report the unstructured condition number $\cond_{|\Omega_{\bsB_1}|}$, the structured effective condition number $\cond_{\rm eff} (\bsA(\Omega_{QS}),B_1(\Omega_{\bsB_1}))$, and structured condition number $\cond^{QS}_{|\Omega_{\bsB_1}|}$ in quasiseparable representation for the multiple right-hand side linear system $\bsA\bsX=\bsB_1$. Furthermore, the unstructured and structured condition numbers $\cond_{|\bsB_2|}$, $\cond_{\rm eff} (\bsA(\Omega_{QS}),\bsB_2)$ and $\cond^{QS}_{|\bsB_2|}$ of the multiple right-hand side linear system $\bsA \bsX= \bsB_2$ are also reported in Table \ref{ta:nx0}. From  Table \ref{ta:nx0}, the  unstructured condition number can be much larger than the structured effective condition number and the structured condition number with respect to the general quasiseparable representation while there are little differences  between the structured effective condition number and the structured condition number with respect to the general quasiseparable representation, which coincide with the collusions of Proposition \ref{pro:the comparison-cond_f and cond_A,f} and Theorem \ref{th:the inequ-cond{QS} and condeff{QS}}. }
\begin{table}
\centering
\caption{\label{ta:nx0} Comparisons of unstructured number and effective structured condition number for multiple right-hand side linear system \eqref{eq:mul} with fixed $n=5$ and $m=2$. }
\label{T1}
{\small
\begin{tabular}{cccccc}\hline
$\cond_{|\Omega_{\bsB_1}|}$& $ \cond_{\rm eff} (\bsA(\Omega_{QS}),\bsB_1(\Omega_{\bsB_1}))$ &$\cond^{QS}_{|\Omega_{\bsB_1}|}$ \\
\hline
 $6.1526\cdot 10^5$&$5.9573 \cdot 10$ & $6.8460\cdot 10$\\
\hline
$\cond_{|\bsB_2|}$& $ \cond_{\rm eff} (\bsA(\Omega_{QS}),\bsB_2)$ &$\cond^{QS}_{|\bsB_2|}$ \\
\hline
 $7.9145\cdot 10^4$&$1.3532\cdot 10 $ & $1.4920\cdot 10$\\
\hline
\end{tabular}}
\end{table}
\end{example}

\begin{example}\label{example:2}

\red{In this example, we test random quasiseparable matrix $\bsA$ and the multiple right hand sides $\bsB$. First we fixed $n=60$ and chose different $m$, where $n$ is the dimensionality of the square matrix $\bsA$ and $m$ is the column numbers of the multiple right hand sides  $\bsB$.  We use the {\sc Matlab}'s built-in function {\sf sprand($n$,$m$,$\rho$)} to generate sparse multiple right hand sides $\bsB_1$ and use the {\sc Matlab}'s command {\sf randn} to generate dense multiple right hand sides $\bsB_2$. The \{1;1\}-quasiseparable $\bsA \in \R^{60\times 60}$ is generated by its   Givens-vector representation via tangent \eqref{eq:GV rep}: $\bsl\in \R^{58}$, $\bsv\in\R^{59}$, $\bsd\in\R^{60}$, $\bse\in \R^{59}$, and $\bsu\in\R^{58}$, where each vectors are obtained by using the {\sc Matlab}'s command {\sf randn}. From Table \ref{ta:nx1}, for different choices of $m$, the  relationship between  $\cond^{GV}_{|\Omega_{\bsB_i }|}$, $\cond^{QS}_{|\Omega_{\bsB_i }|}$, $\cond_{\rm eff}(\bsA(\Omega_{QS}),\bsB(\Omega_{\bsB_i}))$ and $\cond_{|\Omega_{\bsB_i}|}$ are consistent with Theorem \ref{th:the inequ-cond{QS} and condeff{QS}}, Theorem \ref{th:cond(GV)and cond(QS)}, and Theorem \ref{th:cond{GV}-cond{QS}-cond_eff}.}
\begin{table}
\centering
\caption{\label{ta:nx1} Comparison of the unstructured and structured condition numbers of multiple right-hand side linear system \eqref{eq:mul} with a fixed $n=60$ and different choices of $m$.}
\label{T1}
{\small
\begin{tabular}{ccccccc}\hline
 $m$&$\rho$&$\cond^{GV}_{|\Omega_{\bsB_1 }|}$&$\cond^{QS}_{|\Omega_{\bsB_1 }|}$&$\cond_{\rm eff}(\bsA(\Omega_{QS}),\bsB(\Omega_{\bsB_1}))$&$\cond_{|\Omega_{\bsB_1}|}$ \\
\hline
$20$&$0.01$&$1.6207\cdot 10^{2}$&$1.8997\cdot 10^{2}$&$1.5390\cdot 10^2$&$2.0438\cdot 10^2$\\
\hline
$40$&$0.01$&$1.1936\cdot 10^3 $&$1.3766\cdot 10^3$&$1.1429\cdot 10^3$&$1.1258\cdot 10^3$\\
\hline
$60$&$0.01$&$8.5737\cdot 10^1$&$9.3135\cdot 10^1$&$7.2672\cdot 10^1$&$2.7635\cdot 10^2$\\
\hline
$m$& $$&$\cond^{GV}_{|\bsB_2|}$&$\cond^{QS}_{|\bsB_2|}$&$\cond_{\rm eff}(\bsA(\Omega_{QS}),\bsB_2)$&$\cond_{|\bsB_2|}$ \\
\hline
$20$&$$&$2.7182 \cdot 10^2$&$4.0633 \cdot 10^2$ & $3.2123 \cdot 10^2$&$3.1474 \cdot 10^2$\\
\hline
$40$&$$&$9.3369 \cdot 10^2$&$1.3666 \cdot 10^3$&$1.1354 \cdot 10^3$&$1.1166 \cdot  10^3$\\
\hline
$60$&$$&$6.9029 \cdot 10$&$9.7076  \cdot 10$&$7.5664 \cdot 10$&$ 2.8257 \cdot 10^2$\\
\hline
\end{tabular}}
\end{table}
\end{example}

\begin{example}\label{example:3}

In this example, the parameter vectors describing the \{1;1\}-quasiseparable matrix $\bsA$ are generated randomly as follows:
\begin{align}
\bsp  \in \R^{n-1},\, \bsa \in \R^{n-2},\,  \bsq \in \R^{n-1},\,  \bsd \in  \R^{n},\, \bsg  \in \R^{n-1}, \bsb \in \R^{n-2}, \mbox{ and } \bsh \in \R^{n-1}. \label{eq:num para}
\end{align}
 From the numerical results from the general generated random parameters \eqref{eq:num para}, there are marginal differences between the structured effective componentwise condition number and the corresponding unstructured ones. Therefore, in order to make the differences between  the unstructured componentwise condition number \red{$\cond_{|\Omega_{\bsB}|}$ given in \eqref{notaion:3}} and the structured one\,$\cond_{\rm eff}(\bsA(\Omega_{QS}),\bsB(\Omega_{\bsB }))$ to be large, we use the following procedure to construct the parameter vectors \eqref{eq:num para}. First, we randomly select $\ell_1 $ indexes of $\bsp$, and $\ell_2 $ indexes of $\bsa$ are chosen, where $\ell_1 =\lfloor 30\%  \times (n-1)\rfloor$ and $\ell_2 =\lfloor 30\% \times (n-2)\rfloor$. For the selected indexes of $\bsp$ in the ascending order, we multiply the corresponding component of $\bsp $ by the weight $10^{\alpha_i+3}$, where $\alpha_i=1+(i-1)\cdot 4 / (\ell_1 -1)$. On the other hand, in a similar way, we multiply the select component of $\bsa$ by the weight $10^{\beta_i+3}$, where  $\beta_i =\alpha_{\ell_2- i+1}$. The vectors of $\bsd$ and $\bsg $ are rescaled by factors $10^{-3}$ and $10^3$, respectively. We use the command {\tt randn} of {\sc Matlab} to construct vectors $\bsb,\, \bsh$ and $\bsq $, respectively. Under this situation, entries of $\bsb,\, \bsh$ and $\bsq $ satisfy the standard Gaussian distribution and are independent with other entries.

In Table \ref{T1}, let the multi right-hands $\bsB \in \R^{n\times m}$ be generated by using {\sc Matlab}'s command {\tt randn}, which means that $\bsB$ is a dense and unstructured matrix. From Table \ref{T1}, we can conclude that the unstructured condition number \red{$\cond_{|\bsB|}$} given in \eqref{notaion:3}  indeed much larger than the structured one \red{$\cond_{\rm eff}(\bsA(\Omega_{QS}),\bsB)$ defined in Definition \ref{def:cond_eff} when $\Omega_{\bsB}$=$|\bsB|$}, and there exist such multiple right-hand side linear systems with \{1;1\}-quasiseparable coefficient matrices which have ill condition number with respect to perturbations of the entries of the matrix, but have well one with respect to perturbations on the quasiseparable parameters representing the matrix.

Next, we generate the sparse multiple right-hand sides $\bsB$ of \eqref{eq:mul} by using {\sc Matlab}'s built-in function {\sf sprand($n$,$m$,$\rho$)}. Thus $\bsB$ is a random, $n$-by-$m$, sparse matrix with approximately $\rho mn$ uniformly distributed
nonzero entries. The numerical values of \red{$\cond_{|\Omega_{\bsB}|}$ given in \eqref{notaion:1}} and \red{$\cond_{|\Omega_{\bsB }|}^{QS}$}, where $\Omega_{\bsB}$ is a vector composed by the nonzero absolute values of entries of $\bsB$, are reported in Table \ref{T2}. Similar to the observation of Table \ref{T1}, the structured effective condition number \red{$\cond_{\rm eff}(\bsA(\Omega_{QS}),\bsB(\Omega_{\bsB}))$} can be much smaller than the corresponding unstructured condition number \red{$\cond_{|\Omega_\bsB|}$}.\,\,Thus it is necessary to develop structure-preserving algorithms for solving \eqref{eq:mul} when its coefficient matrix is a \{1;1\}-quasiseparable matrix because structure-preserving algorithms can reduce the forward error significantly of the solution with respect to the structured componentwise perturbations on the coefficient matrix $\bsA$ and sparse multiple right-hand sides $\bsB$.

\end{example}

\begin{table}\centering
\caption{\label{ta:nx} The  ratios between \red{$\cond_{|\bsB|}$} and $ \cond_{\rm eff} (\bsA(\Omega_{QS}),\bsB )$ for $n=20,\, 40$, and $60$ with different choices of $m$.}
\label{T1}
{\small
\begin{tabular}{cccccc}\hline
$n$&$m$&$\dfrac{\red{\cond_{|\bsB|}}}{ \cond_{\rm eff}(\bsA(\Omega_{QS}),\bsB )    } $ & $\red{\cond_{|\bsB|}}$  & $ \cond_{\rm eff}(\bsA(\Omega_{QS}),\bsB )  $    \\
\hline
$20$& $4$ & $1.6703\cdot 10^{16} $& $ 2.8277\cdot 10^{17}$&$1.6929\cdot 10^1$ \\
&$8$ & $5.2754\cdot 10^{15} $& $1.0712\cdot 10^{17} $& $ 2.0306\cdot 10^1 $\\
&$12$ & $ 1.3456\cdot 10^{16} $& $7.8757\cdot 10^{17}$ & $5.8529\cdot 10^1$ \\
&$16$& $2.1251 \cdot 10^{16}$ & $ 7.5164\cdot 10^{17}$ & $3.5370\cdot 10^1$ \\
\hline
$40$& $10$ & $1.7507\cdot 10^{16} $& $ 1.8809\cdot 10^{17}$&$1.0744\cdot 10^1$ \\
&$20$ & $6.9995\cdot 10^{15} $& $2.0532\cdot 10^{17} $& $ 2.9334\cdot 10^1 $\\
&$30$ & $ 1.3697\cdot 10^{16} $& $6.7985\cdot 10^{18}$ & $4.9634 \cdot 10^2$ \\
&$40$& $2.8216 \cdot 10^{16}$ & $4.7127\cdot 10^{18}$ & $1.6702\cdot 10^2$ \\
\hline
$60$&$30$ & $7.9458\cdot 10^{15} $& $5.1320\cdot 10^{17} $& $ 6.4587\cdot 10^1 $\\
&$40$ & $ 2.0569\cdot 10^{16} $& $1.2443\cdot 10^{18}$ & $6.0492\cdot 10^1$ \\
&$50$& $1.4127 \cdot 10^{16}$ & $1.9905\cdot 10^{18}$ & $1.4090\cdot 10^2$ \\
&$60$ & $9.9333\cdot 10^{15} $& $1.1307\cdot 10^{18} $& $1.1383 \cdot 10^2  $    \\
\hline
\end{tabular}
}
\end{table}

\begin{table}\centering
\caption{\label{ta:nx} The ratios between $\red{\cond_{|\Omega_{\bsB}|}}$ and $ \cond_{\rm eff} (\bsA(\Omega_{QS}),\bsB(\Omega_{\bsB }))$ for $n=20,40$ and $60$ with different choices of $m$ and $\rho $.}
\label{T2}
{\small
\begin{tabular}{ccccccc}\hline
$n$&$m$&$\rho $&$\dfrac{\red{\cond_{|\Omega_{\bsB}|}}}{ \cond_{\rm eff}(\bsA(\Omega_{QS}),\bsB(\Omega_{\bsB }) )  } $ & $\red{\cond_{|\Omega_{\bsB}|}}$  & $ \cond_{\rm eff}(\bsA(\Omega_{QS}),\bsB(\Omega_{\bsB }) )  $    \\
\hline
$20$&  10 & 0.1 & $3.6276\cdot 10^{15}$ &$6.7572\cdot 10^{16}$  & $ 1.8627\cdot 10^1$  \\
& 10 &  0.3 & $9.9256\cdot 10^{15}$   & $3.4805\cdot 10^{17}$  & $3.5066e\cdot 10^1$ \\
& 10 &  0.5 & $1.2218\cdot 10^{16}$   & $7.5646\cdot 10^{16}$  & $6.1912\cdot 10^0$ \\
& 20 &  0.1 & $ 1.1558\cdot 10^{16}$   & $3.9404\cdot 10^{16}$  &$3.4094\cdot 10^0$ \\
& 20 & 0.3 & $8.6307\cdot 10^{15}$  & $2.2244\cdot 10^{17}$  & $2.5773\cdot 10^1$ \\
& 20 &  0.5 & $8.9736\cdot 10^{15}$  & $1.5430\cdot 10^{17}$  & $ 1.7195\cdot 10^1$ \\
\hline
40& 20 &  0.1 &$9.1705\cdot 10^{15}$   & $4.7061\cdot 10^{17}$   & $5.1317\cdot 10^1$  \\
& 20 &  0.3 & $6.3252 \cdot 10^{14}$    & $2.1572\cdot 10^{16}$ & $3.4105\cdot 10^1$ \\
& 20 &  0.5 & $1.8175\cdot 10^{16}$   & $6.8781\cdot 10^{17}$ & $3.7843\cdot 10^1$ \\
& 40 &  0.1 & $7.0644\cdot 10^{14}$   & $8.1673 \cdot 10^{16}$ & $1.1561\cdot 10^2$ \\
& 40 & 0.3 & $1.3172\cdot 10^{16}$   & $1.9833\cdot 10^{17}$   & $ 1.5057\cdot 10^1$\\
& 40 &  0.5  & $1.8856\cdot 10^{16}$   &$3.0811\cdot 10^{18}$ & $1.6340\cdot 10^2$ \\
\hline
60& 40 &  0.1 & $1.2153\cdot 10^{15}$  & $1.5341\cdot 10^{17}$& $1.2623\cdot 10^2$ \\
& 40 &  0.3 & $1.1415\cdot 10^{16}$ & $3.0454\cdot 10^{17}$& $2.6680\cdot 10^1$ \\
& 40 &  0.5 & $2.0185\cdot 10^{16}$  & $6.8578\cdot 10^{18}$& $3.3974\cdot 10^2$ \\
& 50 &  0.1 & $2.0709\cdot 10^{16}$ & $1.5371\cdot 10^{18}$ &$7.4226\cdot 10^1$ \\
& 50 &  0.3 & $4.4566\cdot 10^{15}$& $2.7132\cdot 10^{18}$& $6.0882\cdot 10^2$ \\
& 50 &  0.5 & $1.1465\cdot 10^{16}$  &$1.7395\cdot 10^{18}$ & $1.5173\cdot 10^2$ \\

\hline
\end{tabular}
}
\end{table}

\section{Concluding remarks}\label{sec:co}

We presented  the explicit expressions for structured condition numbers for multiple right-hand sides linear systems with \{1;1\}-quasiseparable coefficient matrix in the quasiseparable and the Givens-vector representations.  Relationships between different condition numbers for multiple right-hand side linear systems with parametrized coefficient matrices were investigated. We proposed the effective structured condition number for multi-right hand linear systems with \{1;1\}-quasiseparable coefficient matrix. From the numerical experiments, there were some situations that  the effective structured condition number can be much smaller than the unstructured ones. \red{In this paper, we only consider the case that the coefficient matrix $\bsA$ is nonsingular for the multiple right-hand side linear system. We believe  that our finding in this paper can be extended to the case that $\bsA$ is singular where we should consider the condition numbers for linear least squares solution to the multiple right hands linear system \cite{Dai2011} or a more general linear matrix equations \cite{Liao}.} \red{As shown in \cite{Palitta2018}, for Sylvester equations \eqref{eq:syl},   when $\bsA$, $\bsB$ and $\bsC$ are quasiseparable, the solution $\bsX$ is numerically quasiseparable. Therefore, our result can be extended to the condition number study for Sylvester equations  involving quasiseparable structure}. We will report our research progresses  on the above topics elsewhere in the future.

\section*{Acknowledgements }
The authors thank Prof. Dopico  and Dr. Pom{\'e}s for sending \textsc{Matlab} codes of \cite{Dopico2015Structured}. \red{The authors thank reviewers for their valuable comments, which leads to the improvement of the presentation of this paper.} This work is partially  supported  by the Fundamental Research Funds for the Central Universities under the \red{grants 2412017FZ007}.

\end{document}